\documentclass[12pt]{amsart}
\usepackage[dvipdfmx]{graphicx}
\usepackage[all]{xy}
\usepackage{amssymb, amsmath}
\usepackage{setspace}
\usepackage{eufrak}
\usepackage{amsfonts}
\usepackage{float}
\usepackage{amsmath}

\usepackage[super]{nth}
\usepackage{url}
\usepackage{bm}

\calclayout

\usepackage{colortbl}
\usepackage{color}

\def\t{\noindent}

\newcommand{\A}{\mathcal{A}}
\newcommand{\K}{\mathcal{K}}

\newcommand{\Gcal}{\mathcal{G}}
\newcommand{\E}{\mathcal{E}}
\newcommand{\F}{\mathcal{F}}
\newcommand{\C}{\mathbb{C}}
\newcommand{\Z}{\mathbb{Z}}

\newcommand{\Q}{\mathbb{Q}}
\newcommand{\R}{\mathbb{R}}

\newcommand{\Ost}{\mathcal{O}}

\newcommand{\Proj}{\mathbb{P}}

\newcommand\jeden {1\hskip-3.5pt1}

\newcommand\rank{\mbox{\rm rank}\,}
\newcommand\coker{\mbox{\rm coker}\,}
\newcommand\codim{ \mbox{\rm codim}\, }
\newcommand\Dual{\mbox{\rm Dual}\, }

\newcommand\Image{\mbox{\rm Im} }

\def\Hom{{\rm Hom}}

\def\t{\noindent}
\def\rd{\partial}
\def\Hilb{{\rm Hilb}}
\def\CHilb{{\rm CHilb}}

\def\disconst{\alpha_{\mbox{\rm \tiny dis}}}
\def\SM{{\rm SM}}

\def\CH{{\rm CH}}
\def\Aut{{\rm Aut}}

\def\Gr{{\rm Gr}}
\def\Res{{\rm Res}}
\def\git {/\hskip-3.5pt /}

\def\II{I\hskip-1.5pt I}

\def\bw{\bm{w}}

\def\bd{\bm{d}}

\def\b0{\bm{0}}

\newtheorem{thm}{\bf Theorem}[section]
\newtheorem{"thm"}[thm]{\bf `Theorem'}
\newtheorem{cor}[thm]{\bf Corollary} 
\newtheorem{lem}[thm]{\bf Lemma} 
\newtheorem{prop}[thm]{\bf Proposition} 
\newtheorem{definition}[thm]{\bf Definition} 
 
\newtheorem{rem}[thm]{\bf Remark} 
\newtheorem{exam}[thm]{\bf Example}

\def\t{\noindent}

\begin{document}
\setlength{\baselineskip}{15pt}

\title[Thom polynomials]
{Thom polynomials for singularities of maps} 
\author[T.~Ohmoto]{Toru Ohmoto}
\address[T.~Ohmoto]{Department of Pure and Applied Mathematics, Graduate School of Fundamental Science and Engineering, Waseda University,
Tokyo 169-8555, Japan}
\email{toruohmoto@waseda.jp}
%
%
\dedicatory{On the occasion of the \nth{101} anniversary of Ren\'e Thom's birth}
%
%
%
\maketitle
\begin{abstract}
This is a gentle introduction to a general theory of universal polynomials associated to classification of map-germs, called {\em Thom polynomials}. 
The theory was originated by Ren\'e Thom in the 1950s and has since been evolved in various aspects by many authors. In a nutshell, this is about intersection theory on certain moduli spaces, say `classifying spaces of mono/multi-singularities of maps', which  provides consistent and deep insights into both classical and modern enumerative geometry with many potential  applications. 
\end{abstract}
%
%

\setcounter{tocdepth}{3}
\tableofcontents

\section{Introduction}\label{introduction}

{\it Thom polynomials} are a key component of a general enumerative theory for singularities of real and complex mappings -- these are universal cohomological obstructions to the appearance of singular points of prescribed types in given mappings. 
Roughly speaking, a {\em singularity} of a differentiable map $f$ means a point where the differential $df$ is degenerate (and the germ of $f$ at such a point) or also a multiple image point of the map. There are a lot of different local types of singularities of maps, and the locus of those singular points carries important information about global topological nature of the source and the target spaces. Historically, the origins lie in classical enumerative geometry founded by pioneers of the 19th century, and based on it, the theory was born as a branch of singularity theory of differentiable maps in the mid-20th century. In fact, Ren\'e Thom simultaneously developed singularity theory, cobordism theory, characteristic classes and cohomology operations, armed with his transversality theorem, since his doctor thesis \cite{Thom51} -- especially, in 1957, the year before he received the Fields Medal, Thom gave the inaugural lecture for professorship at the University of Strasbourg, in which he manifested there is a deep correspondence from local classification theory of map-germs to global theory on topology of manifolds and maps \cite{Thom57}. This is what people later called Thom polynomial theory\footnote{The author heard this story from B. Teissier at the workshop ``Hommage \`a Ren\'e Thom, IRMA, Strasbourg, 2016.}. 
Nowadays, the theory itself has been expanding to much broader fields dealing with various enumeration problems in algebraic geometry, differential topology, complex analytic geometry, representation theory, and their interactions.

Excellent surveys of the theory up to the early 1990s may be found in Kleiman \cite{Kleiman76, Kleiman77, Kleiman87} from algebraic geometry side, and in the fourth chapter of \cite{AGLV} written by Vassiliev (Vasil'ev) from topology side. The present article aims to briefly outline subsequent developments with emphasizing new ideas --  {\em equivariant cohomology}, {\em  localization} and {\em Hilbert schemes of points} -- and how those are merged with old traditional methods.  As the theory still now continues to develop, this survey may include  conjectural (but mostly reliable) statements and should therefore be considered an interim report at this stage. Additionally, content selection is heavily biased by the author's interests -- in particular, due to space constraints, we restrict our attention only to some basic topics in singularity theory of maps. 
Consequently, there are several important topics which are not picked up or less mentioned in this survey. 
For instance, we will not enter advanced topics in contemporary Schubert calculus and enumerative geometry on various moduli stacks, which are deeply connected to geometric representation theory, quantum integrable systems, mirror symmetry, Gromov-Witten and Donaldson-Thomas theories and so on. 
Also about real Thom polynomials for $C^\infty$-mappings in differential topology,  the range of related topics are quite large, e.g., from Gromov's h-principle to low dimensional topology including smooth Poincar\'e conjecture, so we will only tough on a few basic matters. 
For these topics not covered in this survey, readers are referred to appropriate references. 

A new excellent survey paper on Thom polynomial theory written by R. Rim\'anyi \cite{Rimanyi24} recently appeared. It differs in style, focus and perspective from this survey and is highly recommended reading.  Readers are also referred to Rim\'anyi's website  \cite{RimanyiPortal}, which contains many latest computational results of Thom poylnomials. 

We work mainly with the complex analytic category (and real $C^\infty$ category in the last section). It is also possible to take algebraic geometry context with some appropriate modifications. We assume that readers are familiar with basics of differential topology or algebraic geometry, and especially, preliminary elements for the following subjects (e.g., see the very beginning chapters in cited survey papers included in Handbook of Geometry and Topology on Singularities III): 
\begin{itemize}
\item 
Thom-Mather theory of singularities of maps (cf.  \cite{MondNuno22}); 
\item 
characteristic classes of vector bundles  (cf.  \cite{Brasselet22}); 
\item 
intersection theory in algebraic geometry (cf.  \cite{Al22}). 
\end{itemize}

 Throughout this paper,  $H_*$ stands for the Borel-Moore homology and $H^*$ the singular cohomology, where the coefficients are integers unless specifically noted. 
For complex manifolds,  $H_*$ and $H^*$ are identified via the Poincar\'e duality as usual. 
In algebraic geometry context, we use Chow groups $\CH_*$ (Chow rings $\CH^*$ for non-singular varieties) as in  \cite{Fulton}. 

The {\em codimension $\kappa$ of a map $f: M \to N$} between manifolds is defined as $\kappa:=\dim N -\dim M$ 
and also the relative dimension to be $d:=-\kappa$. 

For a complex vector bundle $E \to M$, the Chern class $c_i(E) \in H^{2i}(M)$ is assigned ($c_0(E)=1$, $c_i(E)=0$ for $i > \rank E$). 
We often use the {\em quotient Chern class} $c_i=c_i(f)$ associated to a map $f: M \to N$, which is defined  
as the $i$-th Chern class of the difference bundle $f^*TN-TM$, i.e., 
$$c(f):=1+c_1(f)+c_2(f)+\cdots =\frac{1+f^*c_1(TN)+f^*c_2(TN)+\cdots}{1+c_1(TM)+c_2(TM)+\cdots} \;\; \in H^*(M).$$ 
We also use $\bar{c}_i(f):=c_i(TM-f^*TN)$, i.e., 
$$\bar{c}(f):=1+\bar{c}_1(f)+\bar{c}_2(f)+\cdots \;\; \mbox{with} \;\; c(f)\bar{c}(f)=1.$$
Note that every $\bar{c}_i(f)$ is a polynomial in $c_1(f), \cdots, c_i(f)$,  vice versa, and the dual version is $c_i(T^*M-f^*T^*N)=(-1)^i \bar{c}_i(f)$. 

If $f$ is proper, the pushforward $f_*: H^*(M) \to H^{*+\kappa}(N)$ is defined (precisely, it is the Gysin homomorphism through the Poincar\'e duality; we simply denote it by $f_*$ unless any confusion occurs; the superscript should be multiplied by $2$ in complex case). 
It satisfies the projection formula  
$$f_*(\alpha \cdot f^*\beta)=f_*(\alpha) \cdot \beta.$$ 
Since $f_*$ is just an additive homomorphism, the {\em Landweber-Novikov class} with index $I=(i_1, i_2, \cdots, i_k)$  
$$s_I(f):=f_*(c^I(f))=f_*(c_1(f)^{i_1}c_2(f)^{i_2}\cdots c_k(f)^{i_k})$$
takes important role, especially in complex cobordism theory. 
We often use the same notation $s_I$ to mean its pullback $f^*s_I$ (i.e., drop the letter $f^*$), for it is easily understood from the context. Also the readers should not confuse the notation with Segre classes of (virtual) vector bundles.

\subsection{Two guiding examples} 
As an initial guide to readers, we first present typical examples of Thom polynomials for classifications of linear maps and for non-linear map-germs. 
Throughout this survey, these examples will be referred to repeatedly.

\subsubsection{Thom-Porteous formula} \label{thom_porteous}
Consider the classification of linear maps $h: \C^m \to \C^n$ under linear coordinate changes of the source and the target spaces. The orbits are determined by the rank, or equivalently, the dimension of kernel. Note that using $\dim \ker h$ is better than $\rank h$ in the sense that $\dim \ker h$ does not change under the trivial suspension $h: \C^m \to \C^n \leadsto h \times id_k:  \C^{m+k} \to \C^{n+k}$. 
The corresponding Thom polynomials are described as follows. 

Let $E \to M$ and $F \to M$ be complex vector bundles 
of rank $m$ and $n$, respectively, over the same base manifold $M$. 
Let $\varphi: E \to F$ be 
a vector bundle morphism, that is a section of the vector bundle $\Hom(E,F)=E^*\otimes F \to M$. 
We are herewith interested in the loci corresponding to singular linear maps; put 
$$\Sigma^k(\varphi)=\{\; x \in M\; | \; \dim \ker \varphi_x  =  k\; \}$$
and define the {\it degeneracy locus} of $\varphi$ to be the closure 
$$\overline{\Sigma^k(\varphi)}=\bigcup_{s\ge k} \Sigma^s(\varphi).$$ 
The expected codimension of the locus is given by $kl$ where $k=\dim \ker$ and $l:=\dim \coker=n-m+k\, (=\kappa+k)$. 
In particular, for a map $f: M^m \to N^n$, 
we denote by $\Sigma^k(f)$ the locus for the differential $df: TM \to f^*TN$,  called the {\em first order Thom-Boardman singularity stratum}. 

The {\it Giambelli-Thom-Porteous formula} \cite{Porteous, Fulton} states that for suitably generic $\varphi$, 
the dual to the degeneracy locus is expressed 
 by 
\begin{equation}\label{GTP}
{\rm Dual}\, [\overline{\Sigma^k(\varphi)}]  
 \; = \; \Delta_{l^k}(c(F-E)) \; = \; 
\overbrace{
\left|
\begin{array}{ccc}
c_{l} & c_{l+1} & \cdots \\
c_{l-1} & c_{l} & \cdots \\
\vdots & \vdots & \ddots \\
\end{array}
\right|
}^k
\end{equation}
where $c_i$ is the $i$-th component of the quotient Chern class $c(F-E)$. 
The Schur function of type $l^k=(l, \cdots, l)$ in the right hand side of (\ref{GTP}) 
is the {\em Thom polynomial $tp(\Sigma^k)$ for the type $\Sigma^k$}. 
Even for non-generic case, i.e., for an arbitrary bundle map $\varphi$, the formula is valid after replacing the left hand side of (\ref{GTP}) by a certain intersection class, called the {\em degeneracy loci class}, see Fulton \cite{Fulton}. 
We remark that there is a well-known identity via conjugacy: 
$$\Delta_{\lambda}(c(F-E)) = \Delta_{\lambda^\vee}(c(E^*-F^*))$$
(cf. \cite[Lem.14.5.1]{Fulton}; note that the notation is inconsistent across different references). 

In case that $\kappa \ge 0$ and $E$ is trivial, then the largest degeneracy locus represents $c_{\kappa+1}(F)$, i.e., 
the Chern class of a vector bundle itself is a sort of Thom polynomial. 

There are several directions for generalization of the Thom-Porteous formula. 
When replacing the group by Borel subgroups, we meet the (matrix) Schubert calculus, so this passage goes to geometric representation theory and algebraic combinatorics (\S \ref{thom_porteous2} and \S \ref{modern_Schubert}). 
When one considers non-linear mappings, it goes to singularity theory of maps, see the next example.

\subsubsection{Fold and Cusp}  \label{fold_and_cusp} 
Recall the {\it classical Riemann-Hurwitz formula}. Let $f:M \to N$ be a surjective holomorphic map between compact nonsingular complex curves: to each point $x \in M$, we assign $\mu=\mu(f,x)$ so that $f$ is written by $x \mapsto x^{\mu+1}$ in local coordinates centered at $x$ and $f(x)$. Then it holds that 
\begin{eqnarray*}
\int_M \mu(f)  &=&\deg f\cdot \chi(N) - \chi(M)\\
&=& c_1(TN) \frown f_*[M] - c_1(TM) \frown [M]\\  
&=&  c_1(f^*TN-TM) \frown [M] \\
&=& c_1(f) \frown [M].
\end{eqnarray*}
Here we use the fact that the Euler characteristic $\chi(M)$ is the degree of the top Chern class $c_1(TM)$ of the tangent bundle. 
If $f$ is generic, i.e., it has only simple branches ($\mu=1$), called the {\em $A_1$-singularity}, then the formula counts the number of the singular points. 
For a holomorphic map $f: M \to N$ between complex surfaces, if it is appropriately generic, local singularities of $f$ are classified into two types,  called {\em fold} and {\em cusp}, up to local coordinate changes of $M$ and $N$ at the point, 
whose normal forms are respectively given by 
$$A_1: (x, y) \mapsto (x^2, y), \quad A_2: (x, y) \mapsto (x^3+yx, y).$$
In general, for map-germs $(\C^m, 0) \to (\C^{m+\kappa}, 0)$ with $\kappa \ge 0$, the symbol $A_\mu$ denotes the singularity type whose local algebra  is isomorphic to $\C[x]/\langle x^{\mu+1}\rangle$, see \S \ref{classification}. 
 
 \begin{figure}
\includegraphics[height=4.5cm, pagebox=cropbox]{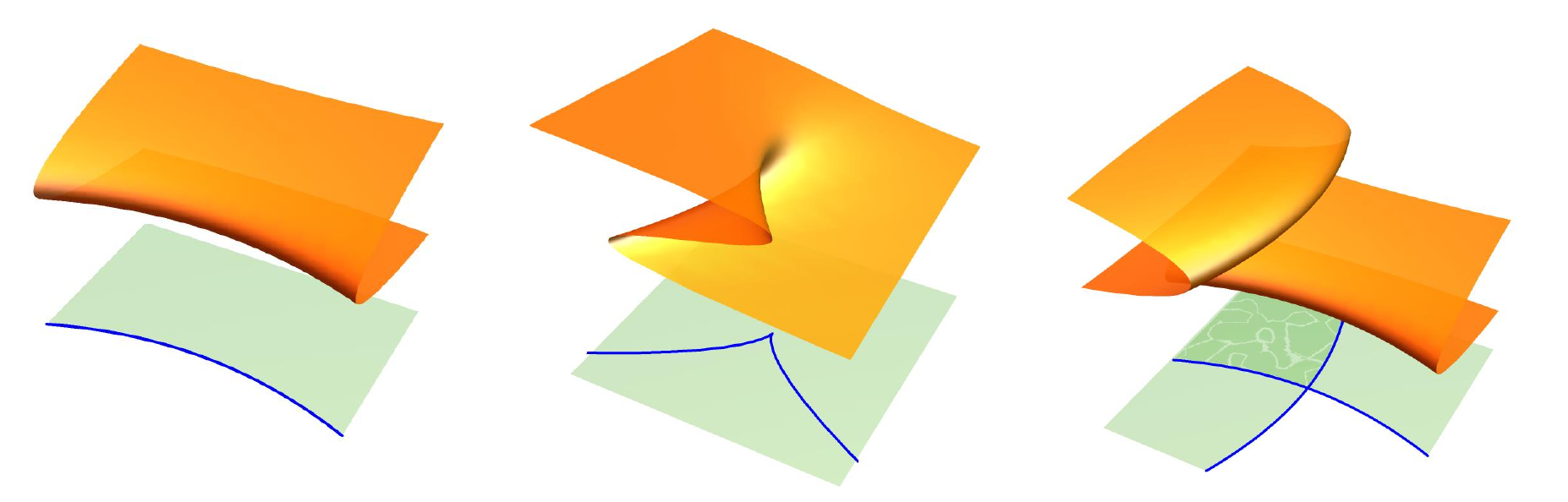}
\caption{Singularities of plane-to-plane maps: fold, cusp and double-folds}
 \end{figure}
 
The critical locus of $f$ is the closure $\overline{A_1(f)}$, which is formed by fold points and isolated cusp points. 
 In his paper  \cite{Thom55}, R. Thom showed (precisely saying, in real case, using the Stiefel-Whitey classes) that the cohomology class represented by the critical locus and the number of cusps are expressed by 
\begin{equation}\label{thom1}
\Dual  [\overline{A_1(f)}] = c_{1}(f), \quad \Dual  [\overline{A_2(f)}] = c_{1}(f)^2+c_{2}(f)
\end{equation}
in $H^*(M)$. The former one is essentially the same as the above classical Riemann-Hurwitz formula. 
In fact, the formulas (\ref{thom1}) hold for any generic maps between manifolds of the same dimension $\ge 2$.  
The right hand side, $tp(A_1):=c_1$ and $tp(A_2):=c_1^2+c_2$, are respectively called {\em Thom polynomials for $A_1$ and $A_2$-singularity type with $\kappa=0$}. The shape of the polynomial depends only on the singularity type, and the information of the map $f$ itself is only encoded in substituting $c_i(f)$ to variables $c_i$ of the polynomial. 
In case of $\kappa=0$, more samples are: 
\begin{align}\label{thom2}
tp(A_3)&= {c}_1^3+3{c}_1{c}_2+{c}_3, \notag \\
tp(A_4)&= {c}_1^4+6{c}_1^2{c}_2+2{c}_2^2+9{c}_1{c}_3+6{c}_4,\\
tp(A_5)&= c_1^5 + 10c_1^3c_2 + 25c_1^2c_3 + 10c_1c_2^2 + 38c_1c_4 + 12c_2c_3 + 24c_5 \notag\\
tp(A_6)&= c_1^6+15c_1^4c_2+55c_1^3c_3+30c_1^2c_2^2+141c_1^2c_4+79c_1c_2c_3 \notag \\
&\qquad \quad +5c_2^3+202c_1c_2+55c_2c_4+17c_3^2+120c_6 \notag \\
& \cdots  \;\; \cdots, \notag
\end{align}
and for singularity types of corank greater than one, e.g.,  $I_{a, b}: (x^a+y^b, xy)$,  
$$tp(I_{2,2})=c_2^2-c_!c_3, \quad 
tp(I_{2,3})=2c_1c_2^2-2c_1^2c_3 +2c_2c_3-2c_1c_4,$$ 
and so on -- 
a full list up to codimension $\le 8$ with $\kappa=0$ can be seen in Rim\'anyi \cite{Rimanyi01}. 
As a catalog of various kinds of Thom polynomials, many computational results have been collected in Rim\'anyi's website `Thom polynomial Portal'  \cite{RimanyiPortal}.

We may find general formulas which contains $\kappa \ge 0$ as a parameter. 
The simplest Thom-Porteous formula says that the Thom polynomial for $A_1$-singularity (i.e., dual to the critical locus) is given by \begin{equation}\label{thom3}
tp(A_1)=\Delta_{(\kappa+1)^1}=c_{\kappa+1}
\end{equation}
and  by a result due to Ronga \cite{Ronga72} we know 
\begin{equation}\label{thom4}
tp(A_2)=c_{\kappa+1}^2+\sum_{i\ge 1}2^{i-1}c_{\kappa+1-i}c_{\kappa+1+i}
\end{equation}
where we read $c_j=0$ for $j<0$. 
Those lead to a general theory of  {\em Thom series}, see \S \ref{tp_series}. 
In particular, from  (\ref{thom1}), (\ref{thom2}), (\ref{thom3}), (\ref{thom4}) and further computations, Rim\'anyi has conjectured the {\em Chern basis positivity for $A_\mu$-singularity types} ($\kappa \ge 0$), that means that all coefficients of Chern monomials appearing in $tp(A_\mu)$ are non-negative \cite{Rimanyi01}. That is still open, and sounds mysterious, as the conjecture is unexpectedly related to complex hyperbolic geometry, see \S \ref{hyperbolic}.   

Let us look back again at a generic map $f: M \to N$ between surfaces. 
There is one more stable singularity type, which is a bi-singularity of fold types, denoted by $A_1^2=A_1A_1$ and  called {\em transverse double-folds}, 
which corresponds to a {\em node} of  the discriminant curve $D(f):=f(\overline{A_1(f)})$ in the target space $N$. 
Assume that $f$ is proper, then the number of nodes of $D(f)$ is half the number counted by 
\begin{equation}\label{thom5}
\Dual [\overline{A_1^2(f)}]=-4c_1(f)^2-2c_2(f)+ c_1(f)f^*s_1(f),
\end{equation}
see \S \ref{multising}.  
Furthermore, the Euler characteristic of $D(f)$ is also universally expressed by 
\begin{equation}\label{thom6}
\chi(D(f))=\int_N ( c_1(TN) s_1(f) + s_{01}(f) - \frac{1}{2} s_1(f)^2),
\end{equation}
see \S \ref{higher_multising}. 
Namely,  the way of the appearance of singularities of prescribed type of maps has its own universal rule, that is dominant by local symmetry of the singularity type. Those kinds of enumerative problems are treated in the Thom polynomial theory, and indeed appear in various fields of applications.  

Finally we mention old examples about fold and cusps of maps between surfaces \cite{Salmon, SempleRoth} --  perhaps, R. Thom must have been familiar with these treasures of French classical geometry since Monge and Poncelet. First, consider a smooth plane curve $C$ of degree $d \ge 2$. Let $M$ be the point-line incidence surface (i.e., it consists of pairs $(p, l) \in C \times \Proj^{2*}$ with $p \in l$) and $f: M \to N:=\Proj^{2*}$ the projection, then the discriminant $D=D(f)$ is just the projective dual to $C$. 
A flex (inflection point) of $C$ corresponds to a cusp of $D$ and a double tangent line to $C$ corresponds to a node of $D$. Applying (\ref{thom1}) and (\ref{thom5}), we recover the Pl\"ucker formula:  
$$\deg D=d(d-1), \;\; \# {\rm Cusps}= 3d(d-2), \;\; \# {\rm Nodes}=d^4-2d^3-9d^2+18d.$$ 
Second, suppose that we are viewing a surface $M \subset \Proj^3$ of degree $d$ from a general viewpoint. It is a map projecting $M$ to a screen $N:=\Proj^2$ and $D(f)$ is the apparent contour, and then a handy calculation of (\ref{thom1}) and (\ref{thom5}) immediately recovers old formulas:  
$$\deg D =d(d-1), \; \# {\rm Cusps} =d(d-1)(d-2), \; \# {\rm Nodes} =\frac{1}{2}d(d-1)(d-2)(d-3).$$ 
Like as those, there are many different geometric situations where fold, cusps, double-folds appear naturally, while our formulas enumerate them in exactly the same way -- namely, the Thom polynomial unites many such counting problems into a form that is relevant only for the type of local singularity. Furthermore, this theory reveals the intrinsic meaning of classical numerical characters of projective varieties. Although it is little-known, Semple-Roth \cite{SempleRoth} wrote that for a projective surface there are exactly four basic characters which generate all other characters of Salmon-Cayley-Zeuthen, and afterwards, Roth claimed an analogous statement for a projective $3$-fold that there are seven basic characters which recover dozens of others. Through computing Thom polynomials \cite{Sasajima17}, it turns out that those basic characters correspond to Landweber-Novikov classes in modern language;  there are four $s_0, s_1, s_2, s_{01}$ for a surface, and seven $s_0, s_1, s_2, s_3, s_{01}, s_{11}, s_{001}$ for a $3$-fold, indeed. This suggests cobordism theory behind (see \S \ref{multising_tp}).  A more advanced application to enumerative geometry, as well as unexpected applications to other fields, will be reviewed in \S \ref{applications}.

\section{Thom polynomials for singularities of maps}\label{sec1}

\subsection{Classification of map-germs}\label{classification}
To begin with, we describe a few basic notions in the Thom-Mather theory, see e.g., \cite{MondNuno22}, \cite[Part I]{MondNuno}, \cite[Part I]{AGV} for details. See also J. Mather's \cite{Mather} for algebraic geometry setting. 

Let $\Ost_m$ be the local ring of holomorphic function germs $(\C^m,0) \to \C$ 
with the maximal ideal $\mathfrak{m}_m=\{h \in \Ost_m, f(0)=0\}$. 
Put $\E(m,n)$ to be the $\Ost_m$-module of all homolorphic map-germs 
$f: (\C^m,0) \to \C^n$, and also put $\E_0(m,n):=\mathfrak{m}_m \E(m,n)$ which consists of all germs $f: (\C^m, 0) \to (\C^n,0)$. 

\

\t{\bf Equivalence}: 
The group of biholomorphic germs $(\C^m, 0) \to (\C^m,0)$ 
is denoted by $\Aut (\C^m,0)$. 
There are two different kinds of natural equivalence relations on map-germs: 
for $\mathcal{G}=\A, \K$, a $\mathcal{G}$-equivalence class of map-germs is called a {\em $\mathcal{G}$-singularity type} ($\mathcal{G}$-type, for short). 

\begin{itemize}
\item {\bf $\A$-classification} (right-left equivalence) classifies map-germs up to isomorphisms of the source and the target spaces. The {\it right-left group} $\A \; (=\A_{m,n})$ is the direct product $\Aut(\C^m,0)\times \Aut(\C^n,0)$, 
which acts on $\E_0(m,n)$ by 
$$(\sigma, \tau).f:=\tau \circ f \circ \sigma^{-1}.$$
Two map-germs $f, g: (\C^m, 0) \to (\C^n,0)$ are {\em $\A$-equivalent} if $f$ and $g$ belong to the same $\A$-orbit. 
\item {\bf $\K$-classification} (contact equivalence)  
classifies zero loci $f^{-1}(0)$ as schemes up to isomorphisms, i.e., it classifies the local algebra 
$$Q(f):=\Ost_m/\langle f_1, \cdots, f_n\rangle_{\Ost_m}$$
where $f_i$ are component functions of $f$. 
In other words, the $\K$-equivalence measures the {\em  tangency of the graphs $y=f(x)$ and $y=0$ in $\C^m \times \C^n$}; the {\it contact  group} $\K \; (=\K_{m,n})$ consists of pairs $(\sigma, \Phi)$ of $\sigma \in \Aut(\C^m,0)$ and map-germs $\Phi: (\C^m,0) \to GL(n,\C)$, and it acts on $\E_0(m,n)$ by 
$$((\sigma, \Phi).f)(x)=\Phi(x)f(\sigma(x)).$$
Two map-germs $f, g: (\C^m, 0) \to (\C^n,0)$ are {\em $\K$-equivalent} if $Q(f) = \sigma^*Q(g)$ for some $\sigma \in \Aut(\C^m,0)$, also equivalently, $f$ and $g$ belong to the same $\K$-orbit. 
In particular, if $f\sim_\A g$, then $f \sim_\K g$, i.e.,  
$\A.f \subset \K.f$. 
\end{itemize}

\

\t {\bf Deformations and stability}: 
Let $f \in \E_0(m,n)$. 
An {\it infinitesimal deformation} of $f$ is a vector field-germ along $f$. 
The space of infinitesimal deformations 
 is regarded as the `tangent space' of $\E(m,n)$ at $f$, and 
is denoted by 
$$\textstyle \theta(f) = \Ost_m \otimes \C^n.$$
For the identity map $id_m$ of $\C^m$, 
$\theta_m:=\theta(id_m)$ is the space of germs of vector fields on $\C^m$ at the origin.  
We set $tf:\theta_m \to \theta (f)$  by $tf(v):=df(v)$ and $\omega f:\theta_n \to \theta(f)$ by $\omega(w):=w\circ f$.  
For our convenience, we put 
$$T\A_e.f=tf(\theta(id_m))+\omega f(\theta(id_n)), \quad  
T\K_e.f=tf(\theta(id_m))+f^*\mathfrak{m}_m\theta(f),$$ 
which are  subspaces of $\theta(f)$ consisting of deformations caused by the actions of source and target isomorphisms 
and contact group action (not necessarily preserving the origin). 
Indeed, the subscript $e$ indicates the extended action of reparametrization groups  
possibly moving the origin, and  the above subspaces are the tangent spaces of the orbits. 

We say that  $f: (\C^m, 0) \to (\C^n,0)$ is  an ({\it infinitesimally}) {\it stable germ}  if it holds that 
$$\theta(f)=T\A_e.f,$$
namely, any deformation of $f$ is re-produced by deformations of coordinate changes of the source and the target spaces. 
It is known, as one of major theorems of J. Mather,  that  for stable germs $f$ and $g$, $f \sim_\K g$ if and only if  $f \sim_\A g$. 
Namely,  for a stable germ, $\A.f = \{\mbox{\small Stable germs}\} \cap \K.f$. 

\

\t {\bf Determinacy}: 
We denote the $k$-jet of $f: (\C^m, 0) \to (\C^n,0)$ by $j^kf(0)$ and the $k$-jet space by $J^k(m,n)$ (the space of Taylor coefficients of order $\le k$ off the constant terms). 
Let $\Gcal=\A$ or $\K$. 
A map-germ $f$ is {\it $k$-$\Gcal$-determined} if every map-germ $g$ with $j^kg(0)=j^kf(0)$ is $\Gcal$-equivalent to $f$. 
For instance, a stable germ $f: (\C^m, 0) \to (\C^n,0)$ is $(n+1)$-$\A$-determined. 
Also $f$ is {\em $\Gcal$-finitely determined} (or {\em $\Gcal$-finite} for short) if it is $k$-$\Gcal$-determined for some $k$. 
Finite determinacy is equivalent to that the orbit $\Gcal.f$ has finite codimension, that is also equivalent to that {\em $\Gcal_e$-codimension} $\dim_\C \, \theta(f)/T\Gcal_e.f$ is finite. In particular, every stable germ has $\A_e$-codimension zero.

\

\t {\bf Jet space and jet extension}:   
The {\em $k$-jet bundle} $J^k(M, N) \to M \times N$ is the fiber bundle with fiber $J^k(m,n)$ and structure group $\A_{m,n}$. 
We define  the {\em $k$-jet extension map} 
$$j^kf: M \to J^k(M, N)$$ 
by assigning to every point $x \in M$ the pair of 
$(x, f(x))$ with the jet of the germ $f: (M, x) \to (N, f(x))$. 
Let $k$ be high enough ($k \ge n+1$). 
The following properties are equivalent: 
\begin{enumerate}
\item  the germ of $f: M \to N$ at $x$ is stable;
\item $j^kf$ is transverse  to the $\A$-orbit at $x$; 
\item  $j^kf$ is transverse  to the $\K$-orbit at $x$.
\end{enumerate}

\

\t {\bf Singularity locus of prescribed type}: 
Throughout the present paper, we use the notations $\eta, \xi, \tau, \cdots$ to mean 
$\K$-invariant subsets of $J^\infty(m,n)$. 
A particular one is a $\K$-singularity type of mono-germs (also a moduli stratum of $\K$-types). 
We say that a mono-germ $f$ is of type $\eta$ if its jet belongs to $\eta$. 
For a holomorphic map $f: M \to N$ ($\dim M=m$, $\dim N=n$), we set 
$$\eta(f):= \{\; x \in M \; | \; 
\mbox{the germ $f$ at $x$ is of type  $\eta$} \; \}$$
and call the closure $\overline{\eta(f)} \subset M$ the {\it $\eta$-type singularity locus} of $f$. 
By finite determinacy, 
we may denote by $\eta \subset J^k(m,n)$ the corresponding $\K^k$-orbit in a jet space of sufficiently high order $k$. 
Since $\eta$ is also $\A$-invariant, 
there is the sub-bundle of the fiber bundle $J^k(M, N)\to M\times N$ with fiber $\eta$, 
denoted by $\eta(M, N)$: 
$$
\xymatrix{
& &J^k(M, N) \ar[d] &  \;\; \overline{\eta(M,N)} \ar@{_{(}->}[l]\\
 \overline{\eta}(f) \;\; \ar@{^{(}->}[r] & M \ar[ur]^-{j^kf} \ar[r]_-{(id, f)} & M \times N &
}
$$
By the definition,  
$\eta(f)=jf^{-1}(\eta(M, N))$ and its closure is 
$$\overline{\eta(f)}=jf^{-1}(\overline{\eta(M, N)})=\bigcup_{\xi \subset \overline{\eta}} \xi(f)$$
where $\xi$ runs over strata of a $\K$-invariant stratification of $\overline{\eta}$, 
i.e., $\xi=\eta$ or boundary strata in $\overline{\eta}-\eta$.  
This is a higher jet version of the degeneracy locus $\overline{\Sigma^k(f)}$ in \S \ref{thom_porteous}.

We say that $f$ is {\em generic with respect to $\eta$} if the jet extension $j^kf$ is transverse to $\eta(M,N)$ and its intersection with the boundary strata is of lower dimension. In this case, especially, the Zariski closure $\overline{\eta(f)}$ forms a cycle in the source manifold $M$, and we are interested in the cohomology classes represented by the locus.

\subsection{Equivariant cohomology and Poincar\'e dual}\label{equiv_cohom}
The $\mathcal{G}$-classification of finitely determined map-germs is essentially reduced to the problem on finite dimensional spaces, i.e., the action of $\mathcal{G}$-equivalence at the level of jets of finite order. So we first review a general framework of {\em Thom polynomials for Lie group action} \cite{Kaz97b, Kaz01, FR04, Ohmoto06}. 

Let $G$ be a complex algebraic Lie group and $X$ a complex algebraic variety on which $G$ acts algebraically. 
In topology, the Borel construction tells that any family $\varphi$ of $G$-orbits in $X$ parametrized by some variety $M$ (precisely,  the pair of a fiber bundle $E \to M$ with fiber $X$ and group $G$ and a section of the bundle) is regarded as a certain morphism, called a  {\it classifying map} 
$$ \bar\varphi: M \longrightarrow   EG\times_G X$$
where $EG \to BG$ is the universal $G$-principal bundle 
and $EG\times_G X=(EG \times X)/G$ is the total space of 
the associated bundle with fiber $X$ and group $G$ over $BG$. 
Any two classifying maps corresponding to isomorphic families are homotopic. 
The {\em $G$-equivariant cohomology of $X$} 
(in the sense of Borel) is defined by 
$$H^*_G(X):=H^*(EG\times_G X).$$ 
For a $G$-morphism $f: X \to Y$, the {\em pullback} $f_G^*: H^*_G(Y) \to H^*_G(X)$ is defined. 

In algebraic geometry context (for reductive $G$), the universal $G$-principal bundle is constructed as ind-varieties  \cite{Totaro, EdidinGraham}. This mimics the construction in topology. 
Roughly saying, take a Zariski open subset $U$ in al linear representation so that $G$ acts on $U$ freely and the quotient variety $U_G:=U/G$ exists, 
then the inductive limit of the quotient map $U \to U_G$ taken over all representations of $G$ (with respect to inclusions) gives an algebraic counterpart to $EG \to BG$. 
For instance, for  $T=\C^*=\C-\{0\}$, we have $U=\C^N-\{0\} \to \Proj^N=U_T$, and for $G=GL_n$, 
let $U$ be an open set in $\Hom(\C^n, \C^N)$ consisting of injective linear maps, 
then $U_G$ is the Grassmaniann of $n$-planes in $\C^N$; those approximates $EG \to BG$. 
For quasi-projective $X$, the diagonal action of $G$ on $X \times U$ is also free and the mixed quotient $X_U:=(X\times U)/G$ exists so that  the projection $p_U: X_U \to U_G$ is a fiber bundle with fiber $X$ and group $G$. 
Then $H_G^*(X)=\varprojlim H^*(X_U)$. 
On the other hand, there is the dual counterpart, the {\em  $i$-th equivariant homology group} 
$H^G_i(X)$, which is constructed in several ways, e.g., a suitable limit of $H_*(X_U)$ through the above Borel construction. 
It satisfies the following properties: 
\begin{itemize}
\item it is a covariant functor to $H_G^*(pt)$-modules: for a proper $G$-morphism $f:X \to Y$, 
the pushforwad $f^G_*: H^G_*(X) \to H^G_*(Y)$ is defined; 
\item it has a specialization to the non-equivariant homology $H_*(X)$; 
\item there exists the {\em equivariant fundamental class} $[X]_G \in H_{2n}^G(X)$ ($\dim X=n$) which defines the homomorphism $\frown [X]_G: H_G^{i}(X) \to H_{2n-i}^G(X)$ so that it is isomorphic if $X$ is non-singular (called the {\em $G$-equivariant Poincar\'e dual}). 
\end{itemize}
We remark again that in case that $X$ is non-compact, we work with (equivariant) Borel-Moore homology. 
For a closed $G$-subvariety $W$ of non-singular $X$ with the inclusion $\iota: W \hookrightarrow X$, 
we write the dual to the image $\iota_*^G[W]_G \in  H_*^G(X)$ by $\Dual_G [W]_G$ (or simply $[W]_G$ for short). 
For a $G$-equivariant proper map $f: X \to Y$ between non-singular $X, Y$, 
we denote the pushforwad though $\Dual_G$ (the Gysin homomorphism) by the same notation 
$f_*^G: H^*_G(X) \to  H^{*+2\kappa}_G(Y)$. 

In particular, we are interested in the case that $X$ is a vector space. Then $EG\times_G X$ is a vector bundle over $BG$, so its $G$-equivariant cohomology is isomorphic to $H^*_G(pt)=H^*(BG)$.

\begin{definition}\upshape \label{tp_G} 
For a $G$-invariant closed subvariety $W$ in a vector space $V$ with $G$-action, 
we define the {\it Thom polynomial} of $W$ to be the equivariant Poincar\'e dual 
$$tp(W):={\rm Dual}^G[W]_G \in H^*_G(V)= H^*(BG).$$
\end{definition}

This notion has implicitly appeared everywhere with different names, e.g., in the setting of torus action, this is known as Joseh polynomial or multidegree etc. 

Let $V$ and $W$ as in Definition \ref{tp_G}. 
Given a vector bundle $\pi: E \to M$ over a complex manifold $M$ 
with fiber $V$ and structure group $G$,  
let  $W(E) \to M$ be the fiber bundle with the fiber $W$ and group $G$. 
For any section $s: M \to E$, the {\em degeneracy loci class of type $W$ for $s$} is defined by
$$m_W(s):=s^*[W(E)] \in H^*(M).$$ 
We say that a section $s: M \to E$ is {\em generic with respect to $W$} if the map $s$ is transverse to every stratum of  (a Whitney stratification of) $W(E)$. 

\begin{prop} \label{prop_tp_G}
\cite{Kaz97b, Kaz01, FR04, Ohmoto06}. 
For a generic section $s: M \to E$ with respect to $\eta$,  it holds that 
$$\Dual\, [W(s)]=m_W(s)=\rho^* tp(W) \in H^*(M)$$
where $W(s):=s^{-1}(W(E))$ and $\rho: M \to BG$ is the classifying map\footnote{Precisely saying, the algebraic classifying map is defined over the total space of an affine bundle over $M$, while this issue does not effect formal computation in cohomology  \cite{Totaro, EdidinGraham}.} 
for  $E \to M$. 
\end{prop}

\subsection{Thom polynomials for $\K$-classification}\label{existence}
$\K$-equivalence of map-germs has a certain stabilization -- essential is the codimension $\kappa=n-m$, rather than $m$ and $n$. 
Two germs $f:(\C^{m+s},0) \to (\C^{n+s},0)$ 
and $g:(\C^m,0) \to (\C^{n},0)$ 
are {\it stably $\K$-equivalent} 
if $f$ is $\K$-equivalent to  
the trivial unfolding $g\times id_s$ with $s$ parameters. 
The following key lemma is easily checked \cite{Damon, Ronga72, Ohmoto94}: 
\begin{lem} \label{stabilization} 
The natural embedding of jet spaces 
$$\Psi: J^k(m,n) \to J^k(m+s, n+s), \quad \Psi(jg(0)):=j(g\times id_{{s}})(0)$$
is transverse to any $\K$-orbits in $ J(m+s,n+s)$. 
 \end{lem}

Apply Definition \ref{tp_G} and Proposition \ref{prop_tp_G} to the setting 
$$G:=J^k\K_{m,n}, \;\; V:=J^k(m,n), \;\; W:=\overline{\eta}.$$
Note that $J^k(m,n)$ is contractible and  
$G=J^k\K_{m,n}$ is homotopic to the $1$-jets $J^1\K_{m,n}=GL_m \times GL_n$, 
and thus $H_G^*(V)$ which is identifed with 
$$H_G^*(pt)=H^*(BGL_m \times BGL_n)
=\Z[a_1, \cdots, a_m, b_1, \cdots, b_n]$$
where $a_i$ and $b_j$ are Chern classes for the source bundle and the target bundle, respectively. 
Hence,  $tp_G(W)$ is a polynomial in $a_i$'s and $b_j$'s. 
The point is that the stabilization property of $\K$-orbits (Lemma \ref{stabilization}) implies that $tp_G(W)$ is {\em supersymmetric}, and hence it is written in terms of quotient variables $c_i$'s where $c=1+c_1+c_2+\cdots =(1+b_1+\cdots)/(1+ a_1+\cdots)$, see e.g., \cite{Kaz03b, Kaz06, FR04} for the detail. 

\begin{thm} \label{tp} 
\cite{Thom55, HaefligerKosinski, Porteous, Ronga72, Damon}. 
For a $\K$-singularity type $\eta$ of  codimension $\kappa$, 
there exists a unique polynomial 
$tp(\eta) \in \Z[c_1, c_2, \cdots]$
so that for any holomorphic map $f: M \to N$ of relative codimension $\kappa$ which is generic with respect to $\eta$, 
 the singularty locus of type $\eta$ is expressed by 
the polynomial evaluated by the quotient Chern class $c_i=c_i(f)$: 
$$\Dual [\overline{\eta(f)}] = tp(\eta)(c(f)) \;\;   \in H^{2\, \mbox{\tiny $\codim \eta$}}(M).$$
\end{thm}

\begin{definition}{\rm 
We call $tp(\eta)$ {\it the Thom polynomial of $\K$-type $\eta$ of mono-germs with codimension $\kappa$}.  
}\end{definition}

There is a counterpart in real $C^\infty$ category using the Stiefel-Whitney class and the Pontrjagin class, which we will pick up in  \S \ref{real_sing} later (Theorem \ref{tp2}).

\begin{rem}\upshape \label{positivity} 
{\bf (Schur Positivity)} 
By the definition, the universal polynomial $tp(\eta)$ is a linear combination of Chern monomials of degree equal to the codimension of $\eta$, while it also has a different expansion. 
Then a certain {\em positivity property} is shown by Pragacz-Weber \cite{PW07, PW08}, see also a survey \cite{Pr16}. 
Given a partition $\lambda=(\lambda_1, \cdots, \lambda_l)$, i.e., $\lambda_1 \ge \lambda_2 \ge \cdots \ge \lambda_l \ge 0$, 
the Schur function with index $\lambda$ is 
$$S_\lambda(E-F):=[s_{\lambda_i+j-i}(E-F)]_{1 \le i, j \le l}$$
where $1+s_1(E-F)+\cdots=\prod (1-\beta)/\prod (1-\alpha)$ and $\alpha$'s and $\beta$'s are the Chern roots of $E$ and $F$, respectively.  
Then we can expand the Thom polynomial with respect to the Schur function basis 
$$tp(\eta)(c(f))=\sum a_\lambda S_\lambda(T^*M -f^*T^*N) \qquad (a_\lambda \in \Z).$$
Since $\overline{\eta(M,N)}$ is a cone bundle,  we directly apply the Fulton-Lazarsfeld theory of numerical positivity of polynomials in the Chern classes of ample vector bundles \cite[\S 12]{Fulton} to prove the desired positivity $a_\lambda \ge 0$. 
\end{rem}

\begin{rem}\upshape \label{cob_tp} 
{\bf (Thom polynomials in other cohomology theory)} 
Instead of using singular cohomology $H^*$ or Chow theory $\CH^*$ to define Thom polynomials for singularities of maps, we may consider to use other cohomology theories. 
For general oriented cohomology theories, note that there is no {\em a priori} defined fundamental class (orientation class) for a singular variety, thus there is no obvious unique way to define an {\em $\eta$-singularity loci class} -- we need some additional data to find a distinguished class supported on the $\eta$-singularity locus. 
For instance, let us consider complex cobordism $MU^*$, or more suitably, algebraic cobordism $\Omega^*$ in the sense of  Levine-Morel \cite{LevineMorel}, and their equivariant version via the Borel construction \cite{HL13, HV11}. 
Pick up a proper $\K$-equivariant birational map onto an orbit-closure, $\pi_\eta: W_\eta \to \overline{\eta} \subset J^k(m,n)$ with $W_\eta$ smooth. Then, for an appropriately generic map $f: M \to N$, we may define a {\em cobordism $\eta$-singularity loci class} in $\Omega^*(M)$ represented by a proper map $W_\eta(f) \to \overline{\eta}(f) \subset M$, where $W_\eta(f)$ is the fiber product of $\pi_\eta$ and the jet extension map $j^kf$ through the above universal construction. The class depends on the choice of $\pi_\eta$, not only the singularity type $\underline{\eta}$. 
Then, it is automatic to see that this class is universally expressed in terms of cobordism Chern classes \cite{LevineMorel, HL13}, i.e., it lives in  
$$\Omega_G^*(pt)=\Omega^*(BGL_m \times BGL_n)=\Omega^*(pt)[[a_1, \cdots, a_m, b_1, \cdots, b_n]].$$
That should be the {\em cobordism Thom polynomial for a desingularization $\pi_\eta$ of  $\eta$-singularity type}, and it may be supersymmetric provided $\pi_\eta$ admits some stabiization property with respect to trivial suspension (see also for multi-singularity case \S \ref{multising}). 
In fact, a typical example is the Thom-Porteous formula in algebraic cobordism $\Omega^*$, which has been studied in  \cite{Hudson12, HPM20} by using Bott-Samelson resolutions of orbit-closures $\overline{\Sigma^k}$ (historically, in real singularity case (cf. \S \ref{real_sing}), already in the 70s, McCrory  \cite{McCrory} considered $\Z_2$-cobordism Thom polynomials for $\overline{\Sigma^k}$ in connection with cohomology operations). 

By the universality of $MU^*$ or $\Omega^*$, a version of Thom polynomial in other oriented cohomology theory (e.g., $K$-theory) can be deduced from this cobordism Thom polynomial. 
Namely, there are some choices for `equivariant fundamental class' in case other than ordinary cohomology theory $H^*$, and of interest is to find a good one among those choices which behaves nicely. A representation theoretic approach to $K$-theoretic Thom polynomials in expansion with respect to Grothendieck polynomials has been studied, especially for $A_2$-singularity type in detail, by Rim\'anyi-Szenes \cite{RS23} through iterated residue integrals technique (cf. \S \ref{tp_series} later). They point out a subtle relation of their $K$-theoretic fundamental class of the orbit-closure (defined using certain characters) with its equivariant resolution, see also Kolokolnikova \cite{Kolo}.  
\end{rem}

\subsection{Desingularization method}\label{desing}
We make a short remark on a traditional approach, called the {\em desingularization method}. 
At the early stages of the theory, the crux of the problem had been recognized as explicitly constructing `resolution of singularities' for the degeneracy loci of prescribed type of maps (indeed, it is natural from the viewpoint of Remark \ref{cob_tp}). 
This approach directly computes Thom polynomials and essentially does not rely on the existence theorem (Theorem \ref{tp}) mentioned above. For instance,  many cancellations occur in computation, and consequently we see that the Thom polynomial is written in terms of quotient Chern classes $c_i=c_i(f)$.

Let $f: M \to N$ be a holomorphic map. Before Mather's classification theory of map-germs appeared, Thom considered the following filtration of (locally closed) submanifolds of $M$.  
First, set 
$$\Sigma^i(f):=\{\; x \in M \; |\; \dim \ker df(x)=i\; \}$$ 
and suppose that it is non-singular (i.e., the transversality condition on $df$ is assumed). Then take the restriction map of $f$ to $\Sigma^i(f)$, 
and set 
$$\Sigma^{i, j}(f):=\{\; x \in \Sigma^i(f) \; |\; \dim \ker d(f|_{\Sigma^i(f)})(x)=j\; \}.$$ 
Assume again it is non-singular, then take the restriction of $f$ to $\Sigma^{i, j}(f)$ and define $\Sigma^{i,j,k}(f)$ and so on. 
This is the so-called {\em Thom-Boardman singularity type $\Sigma^I$} (TB singularity type) with symbol $I=(i \ge j \ge k \ge \cdots)$. 
That is justified by the fact that there are $\K$-invariant algebraic subsets $\Sigma^I$ in $J^r(m,n)$, defined by Jacobian extension (see e.g., \cite{AGV}), so that if the jet-extension map $j^rf$ is transverse to every $\Sigma^I(M, N)$, then it holds that $\Sigma^I(f)=(j^rf)^{-1}(\Sigma^I(M, N))$. For instance, the stratum $\Sigma^{1,1,\cdots, 1,0}$ ($1$ repeats $k$ times) in $J^r(m, n)$ coincides with the $\K$-orbit of $A_k$-singularity type, i.e., the local algebra is isomorphic to $\C[x]/\langle x^{k+1}\rangle$. However, in general, $\Sigma^I$ is not a single $\K$-orbit, and even they do not form a Whitney stratification. 

On the other hand, there is another geometric way to construct some TB singularity loci, by means of {\em intrinsic derivatives} due to Porteous and Ronga. Let us look at a few examples. 
Note that 
$$J^k(M, N)\simeq  \bigoplus_{i=1}^k \Hom({\rm Sym}^iTM, TN)$$
where ${\rm Sym}^i$ means the $i$-th symmetric power and the structure group is reduced to $GL_m \times GL_n$ (thus $J^k(M, N)$ is a vector bundle). 
Let $p_1$ and $p_2$ denote canonical projection to $M$ and $N$, respectively. 
The first intrinsic derivative is a tautological linear map 
$$\alpha: T_x M \to T_y N  \;\; \mbox{at} \;\; \alpha \in J^1(M, N)$$ 
with $x=p_1(\alpha)$ and $y=p_2(\alpha)$, and the second intrinsic derivative is 
$$\beta^*: \ker \alpha \to \Hom(\ker \alpha, \coker \alpha) \;\; \mbox{at} \;\;  (\alpha, \beta) \in J^2(M, N)$$ 
given by ${\rm Sym}^2\ker \alpha \subset {\rm Sym}^2 T_xM \stackrel{\beta}\to T_yN \to \coker \alpha$. Then $\Sigma^i(M, N)$ is defined as the set of $\alpha \in J^1(M, N)$ with $\dim \ker \alpha = i$, and 
$\Sigma^{i,j}(M, N)$ is the set  of $(\alpha, \beta) \in J^2(M, N)$ with $\dim \ker \alpha = i$ and $\dim \ker \beta^*=j$.  

Let $\pi: G_i \to J^1(M,N)$ be the Grassmann bundle of $i$-dimensional subspaces $\lambda \subset T_xM$, and 
 $\xi_i$ the tautological bundle over $G_i$.  
 The vector bundle $\Hom(\xi_i, \pi^*p_2^*TN) \to G_i$ naturally admits a section $\sigma$ which assigns to $(\alpha, \lambda) \in G_i$  the composed linear map $\lambda \hookrightarrow T_xM \stackrel{\alpha}{\to} T_yN$, and $\sigma$ is indeed transverse to the zero section. Hence, the zero locus $Z(\sigma)$ of this section is smooth, and it is birationally mapped to the closure $\overline{{\Sigma}^i(M, N)}$ via $\pi$. Applying the Gysin map of $\pi$, we have  
 $$\pi_*[Z(\sigma)]=\pi_*c_{top}(\Hom(\xi_i, \pi^*p_2^*TN))=[\overline{{\Sigma}^i(M, N)}],$$ 
 and a computation yields the Thom-Porteous formula (\ref{GTP}) mentioned in \S \ref{thom_porteous}  \cite{Porteous, Ronga72, Fulton}. 

Let $\pi':G_{i,j} \to J^2(M, N)$ be the flag bundle of $\tau \subset \lambda \subset T_xM$ with $\dim \lambda=i$ and $\dim \tau=j$ and $G_i$ be defined over $J^2(M, N)$ (via pullback). Over the pullback of $Z(\sigma)$ via the natural projection $G_{i, j}\to G_i$, 
the desingularization of $\overline{{\Sigma}^{i, j}(M, N)}$ is constructed through a bit technical process using  $\beta^*$ \cite{Ronga72}. 
In particular, a computation of the Gysin image for $\Sigma^{1,1}$ leads to (\ref{thom4}) mentioned in \S \ref{fold_and_cusp}. 

In summary,  this traditional approach consists of two steps. The first one is to explicitly construct a tower of some flag bundles 
$$\pi:G \to \cdots \to G_{i, j} \to G_i \to J^k(M, N)$$
 and a vector bundle $E \to G$ with a section $\sigma$ defined by intrinsic derivatives such that the zero locus $Z(\sigma)$ is smooth and gives a desingularization (birational map) $\pi: Z(\sigma) \to \overline{{\Sigma}^I(M, N)}$.  The second is to directly compute the Gysin image $\pi_*(c_{top}(E))$ to get the Thom polynomial. Here, one may use a general formula of the Gysin map for the projection of (partial) flag bundles given by Damon \cite{Damon73} with residue expression. 

However, this strategy fails quickly. Even for order two, the computation in the second step is not easy at all. For order more than two, it is not known in general whether such a nice desingularization in the first step exists or not, and only a few concrete examples are known in equidimensional case ($\kappa=0$); e.g., by finding desingularizations, Thom polynomials are computed for the closures of 
$\Sigma^{1,1,1}$ \cite{Porteous}, $\Sigma^{1,1,1,1}$  \cite{Gaffney} and $\Sigma^{1,1,1,1,1}$ \cite{Turnbull} 
 (those are exactly $tp(A_\mu)$ with $\mu=3,4,5$ in \S \ref{fold_and_cusp}, respectively) and also see \cite{Damon}. 
To go further, we need a different perspective, which we will describe in \S \ref{equivariant}.

\section{Equivariant localization}\label{equivariant}

\subsection{Symmetry of singularities and interpolation}\label{interpolation} 
Around the end of the last century, a novel efficient approach to computing Thom polynomials without using any desingularization was proposed by R. Rim\'anyi \cite{Rimanyi01}, that brought a breakthrough in this field. 

Let $V$ be a complex affine space with action of a complex algebraic Lie group $G$. 
Assume that $V$ consists of finitely many $G$-orbits $\xi$, and also that for every orbit $\xi$, 
 the Euler class of the normal bundle $\nu_\xi$ of $\xi$ in $V$ is non-zero, $c_{top}(\nu_\xi)\not=0$. 
Let $G_\xi$ be the stabilizer subgroup of $\xi$, and 
$\iota_\xi: G_\xi \subset G$ the inclusion. 
Then it turns out that 
$$\bigoplus \iota_\xi^*: H_G^*(V) = H^*(BG){\longrightarrow} \bigoplus_\xi  H^*(BG_\xi) 
$$
is an isomorphism. In particular it holds that 
\begin{equation}\label{REq}
\iota_\xi^*\; tp(\eta)=
\left\{
\begin{array}{cc}
0 & (\xi \not\subset \overline{\eta})  \\
c_{top}(\nu_\eta) & (\xi=\eta)
\end{array}\right.
\end{equation}
The system of these equations (\ref{REq}) has a unique solution modulo torsion \cite{FR04}. 
The first equations in (\ref{REq}) means that $tp(\eta)$ lives in the kernel of the ring homomorphism $H_G^*(V)\to H_G^*(V-\overline{\eta})$, which is isomorphic to the natural image of 
$$H_G^*(V, V-\overline{\eta})\simeq H_*^G(\overline{\eta})\;\;\;\; \mbox{(equivariant Alexander duality)}.$$
The kernel is called  the {\em avoiding ideal associated to $\eta$} in Feh\'er-Rim\'anyi  \cite{FR04} or the {\em ideal of cohomology classes supported on $\overline{\eta}$} in Fulton-Pragacz \cite{FP98}. 
In this ideal, the above second equation  in (\ref{REq}) distinguishes the $G$-fundamental class of $\overline{\eta}$, that is nothing but $tp(\eta)$. 

Also note that it suffices to take the maximal torus of $G_\xi$ for each $\xi$ by splitting principle. 
Once we write $tp(\eta)$ as a homogeneous polynomial in generators of $H^*(BG)$ with unknown coefficients,  (\ref{REq}) gives the system of linear equations in the unknowns, and solving it determines the explicit form of $tp(\eta)$.  Indeed, usually it is overdetermined. 

Today, this is known as the {\em restriction} or {\em interpolation method} due to R. Rim\'anyi. Originally it comes from a slightly different and deeper background, {\em singular Thom-Pontrjagin construction} \cite{Rimanyi96, RimanyiSzucs} -- namely, it works not only for Lie group action but also for {\em Lie groupoid action} or {\em Artin stack with global stratification}, see \S \ref{multising_tp} later. 

Below, we demonstrate some typical computations using this method, according to \cite{Rimanyi01, FR04, Kaz03b}, for objects treated in Guiding Examples \S \ref{thom_porteous} and \S \ref{fold_and_cusp}. 

\subsubsection{Thom-Porteous formula, revisited}\label{thom_porteous2}
We consider the classification of linear maps by linear coordinate changes. Let 
$$G=GL_m\times GL_n, \;\; V={\rm Hom}\, (\C^m, \C^n)$$ 
and $\eta=\Sigma^k$, the $G$-orbit of linear maps with kernel dimension $k$. 
To determine the universal form of $tp(\Sigma^k)$, we see that it is enough to use only the last equation in (\ref{REq}). 
Put ${l}=n-m+k$.

Take a representative $h$ in $\Sigma^k$ in a standard way, then 
the tangent space of the orbit at the point $h$ consists of 
all matrices in the following form: 
$$h=\left[
\begin{array}{c|c}
I_{m-k} & \; O \\
\hline
O & \; O
\end{array}
\right]
\in \Sigma^k, \quad 
T_h\Sigma^k = 
\left\{ \; 
\left[
\begin{array}{c|c}
\;\; *\; \; & \; *\;  \\
\hline
\;\; *\;\;  & \; O\; 
\end{array}
\right]  \;  \right\} \subset T_hV=V.$$
Note that the normal space is isomorphic to $\Hom(\ker h, {\rm coker}\, h)$. 
The stabilizer group of $h \in \Sigma^k$, denoted by $G_k$, 
consists of pairs of square matrices 
$$\left(
\left[
\begin{array}{c|c}
P  & O  \\
\hline
\; * \;   & A
\end{array}
\right] , 
\left[
\begin{array}{c|c}
P &  \; * \: \\
\hline
O  & B
\end{array}
\right]
\right), \quad (A, B, P) \in GL_k\times GL_l \times GL_{m-k}.$$
Let  $\rho_1, \rho_2, \rho_3$ be the 
representations of $G_k$ on $\ker h=\C^k$, ${\rm coker}\, h=\C^l$, ${\rm Im}\, h=\C^{m-k}$, respectively, 
and then the normal bundle $\nu_k$ of the orbit $\Sigma^k$ in $V$ is obtained from the universal bundle $EG_k \to BG_k$ with $\rho_1$ and $\rho_2$ acting on $\Hom(\C^k, \C^l)$. 
Put $c(\rho_1)=\prod_{i=1}^k (1+\alpha_i)$ and  $c(\rho_2)=\prod_{j=1}^l (1+\beta_j)$. 
Then using (\ref{REq}) and the classical resultant,  we have
\begin{eqnarray*}
\iota_k^*\; tp(\Sigma^k)&=&c_{top}(\nu_k)=c_{kl}(\rho_1^*\otimes \rho_2)\\
 &=&\prod (\beta_j-\alpha_i)=\Delta_{l^k}(c(\rho_2-\rho_1))\, \in H^*(BG_k).
 \end{eqnarray*}
 
Denote by $a_i$ and $b_j$
 universal Chern classes for $GL_m$ and $GL_n$, respectively, and then 
 $$H^*(BG)=\Z[a_1, \cdots, a_m, b_1, \cdots, b_n].$$ 
 The representation of $G_k$ on $\C^m$ and $\C^n$ are 
 $\lambda_1=\rho_1\oplus \rho_3$ and $\lambda_2=\rho_2\oplus \rho_3$. 
Thus 
$$\iota_k^*: H^*(BG) \to  H^*(BG_k)$$ 
is determined by  
$c' \mapsto c(\lambda_1)$ and  $c''\mapsto c(\lambda_2)$, 
and especially,  it sends the quotient variables 
$$1+c_1+c_2+\cdots :=\frac{1+b_1+\cdots + b_n}{1+a_1+\cdots + a_m} 
\longmapsto \frac{c(\lambda_2)}{c(\lambda_1)} =1+c_1(\rho_2-\rho_1)+\cdots. $$
Hence $\iota_k^*\Delta_{l^k}(c)= \Delta_{l^k}(c(\rho_2-\rho_1))$. 
Finally it is checked that $\iota_k^*$  is {\it injective} for degree $\le k {l}$, 
thus the consequence is  
$$tp(\Sigma^k)=\Delta_{l^k}(c).$$ 
Moreover, the avoiding ideal associated to $\Sigma^k$ is explicitly determined, see \cite{FP98, FR04}. 

This method works for actions of Borel subgroups as well, that leads to application to (matrix) Schubert calculus and quiver representations \cite{FR03, FR02b, FR04, BFR05, BR07} and also other topics such as algebraic combinatorics in configuration problems, e.g.,   \cite{FNR05, FNR06, FNR09, FNR12, DFR13} and so on. 
Afterwards, this research direction has been widely expanding, that will be touched in \S \ref{modern_Schubert} later. 

For the second order Thom-Boardman singularity type $\Sigma^{i, j}$, the interpolation has also been discussed by Feh\'er-K\"um\"oves \cite{FK06} based on Ronga's disingularization. It will be treated in a more general framework (\S \ref{tp_series}).

\subsubsection{Fold and Cusp, revisited}\label{fold_and_cusp2}
We demonstrate how to compute the Thom polynomial for $A_2$-singularity type (cusp) 
for maps between the same dimensional manifolds, i.e., in case of $\kappa=0$ (see \S \ref{fold_and_cusp}). 
By the existence theorem of Thom polynomials, it has the form
$$tp(A_2)=A c_1^2+Bc_2$$
in quotient Chern classes $c_i=c_i(\mbox{target}-\mbox{source})$.  
Our task is to determine the unknown integer coefficients  $A, B$. 

The normal form of stable map-germ of type $A_2$  is given by a quasi-homogeneous polynomial map 
$$A_2: \C^2 \to \C^2, \;\; (x,y) \to  (x^3+yx, y).$$
In jet space of order $\ge 3$, the jet-extension of this map is transverse to the $\K$-orbit at the origin, so 
this map itself is regarded as a normal slice to the $\K$-orbit of $A_2$-singularity type. 
The maximal torus $T$ of the stabilizer group gives the torus actions on the source and the target  
$$\rho_0=\alpha\oplus \alpha^{\otimes 2}, \quad 
\rho_1=\alpha^{\otimes 3}\oplus \alpha^{\otimes 2}$$
satisfying $A_2 \circ \rho_0=\rho_1 \circ A_2$, where  $\alpha$ is a character $T \to \C^*$.  
Take the dual tautological line bundle 
$\ell = \Ost_{\Proj^N}(1)$ over 
a projective space $\Proj^N$ of large dimension $N \gg 0$
(or the classifying space $BT=\Proj^\infty$ of the torus $T$).  
Define 
two vector bundles of rank $2$ 
$$E_0 \;(=E_0(A_2)) :=\ell \oplus \ell^{\otimes 2}, \qquad 
E_1\;(=E_1(A_2)):= \ell^{\otimes 3} \oplus \ell^{\otimes 2}.$$ 
Since the normal form of $A_2$ is invariant under the torus action,  
we can glue together local patches $id_{U_i} \times A_2: U_i \times \C^2 \to U_i \times \C^2$. 
The resulting map $f_{A_2}:E_0 \to E_1$ is a stable map between the total spaces 
$E_0$ and $E_1$ so that the following diagram commutes and 
the restriction to each fiber  is $\A$-equivalent to the normal form of $A_2$. 
 We call $f_{A_2}$  the {\it universal map for $A_2$}. 
$$\xymatrix{
E_0 \ar[rr]^{f_{A_2}} \ar[dr]_{p_0} & & 
E_1\ar[dl]^{p_1} 
\\
& \Proj^N &
}
$$
The loci $A_2(f_{A_2})$ and  $f(A_2(f_{A_2}))$) are 
just  the zero sections of $E_0$ and  of $E_1$, respectively. 

Put $a=c_1(\ell)$ and then 
 $c(E_0)=(1+a)(1+2a)$ and $c(E_1)=(1+ 3a)(1+2a)$ in $H^*(\Proj^N)=\Z [a]/(a^{N+1})$. 
Note that $H^*(E_0) \simeq H^*(\Proj^N) \simeq H^*(E_1)$ via $p_0^*$ and $p_1^*$. 
Since the $A_2$-locus in the total space $E_0$ is the zero section, we see 
$$\Dual [\overline{A_2}(f_{A_2})]= c_2(p_0^*E_0)=c_2(E_0) = 2a^2. $$
On the other hand, the tangent bundles $TE_0$ and $TE_1$ of the total spaces canonically 
 split as $TE_i=p_i^*(E_i\oplus T\Proj^N)$, thus $f_{A_2}^*TE_1 - TE_0 = p_0^*(E_1-E_0)$. 
Therefore,
the quotient Chern class for $f_{A_2}$ is written as follows: 
\begin{eqnarray*}
c(f_{A_2})&:=&c(f_{A_2}^*TE_1 - TE_0)
= c(E_1-E_0)=\frac{1+3a}{1+a} =1+2a-2a^2+\cdots.
\end{eqnarray*}
Hence, $c_1(f_{A_2})=2a$, $c_2(f_{A_2})=-2a^2$, 
so we have 
$$tp(A_2)(f_{A_2})=A c_1^2+B c_2 = (4A-2B)a^2.$$
By Theorem \ref{tp}, $tp(A_2)(f_{A_2}) = \Dual [\overline{A_2}(f_{A_2})]$, 
thus $(4A-2B)a^2=2a^2$, i.e., $2A-B=1$.

Next we apply $tp$ to the universal map 
of adjacent singularities. Consider  
$$A_1: \C \to \C, \;\; x \mapsto x^2$$
and the associated universal map $f_{A_1}: E_0 \to E_1$, 
where $E_0=E_0(A_1)=\ell$ and $E_1=E_1(A_1)=\ell^{\otimes 2}$ 
 in the same way as above. 
Obviously, the universal map does not have $A_2$-singularity: 
$A_2(f_{A_1})=\emptyset$, 
thus by Theorem \ref{tp} again, we have 
$tp(A_2)(f_{A_1}) = \Dual[\emptyset]= 0$. 
Since $c(f_{A_1})=\frac{1+2a}{1+a}=1+a-a^2+\cdots$, 
we have $A-B=0$. 
 
Finally, solving linear equations in $A, B$ we get $A=B=1$, 
thus we conclude that 
$$tp(A_2)=c_1^2+c_2.$$
Note that instead of working over $\Proj^N$, one can perform the same computation formally as a map from the representation ring $R(T)$ to $H^*_T(pt)=\Z[a]$. 

By this method, Rim\'anyi first computed Thom polynomials for $\K$-orbits in equidimensional case $\kappa=0$ up to codimension $8$ in \cite{Rimanyi01}, that opened a new stage of this theory. Indeed, 
it is known that any stable-germs appearing in the so-called Mather's nice dimensions admit quasi-homogeneous normal forms, and also that for such stable-germs, there are only finitely many nearby $\K$-orbits (i.e., no moduli strata nearby), see e.g., \cite{MondNuno}. 
Thus, once fixing $\kappa$, there are finitely many $\K$-orbits with non-trivial normal Euler classes {\em up to certain codimension} (determined by Mather's nice dimension), so the interpolation method works well for computing Thom polynomials for such stable singularity types, at least theoretically.

A more difficult question is to seek for general formulas of Thom polynomials involving $\kappa$ as a parameter. For $A_2$-singularity with $\kappa \ge 0$ (i.e., $\Sigma^{1,1}$-type), there was Ronga's result as mentioned in Introduction (\S \ref{fold_and_cusp}), and for  $A_3$-singularity ($\kappa \ge 0$), see B\'erczi-Feh\'er-Rim\'anyi \cite{BFR} and Lascoux-Pragacz \cite{LP10}.  Also there are some other partial results due to Pragacz \cite{Pr05, Pr07, Pr08}, \"Ozt\"urk \cite{Oz07},  \"Ozt\"urk-Pragacz \cite{OP12} by the interpolation method with combinatorial techniques of Schur functions. A more general and modern treatment for this problem will be discussed in the next section.

\subsection{Thom series and punctual Hilbert schemes}\label{tp_series}
The next breakthrough came with torus-equivalriant localization technique by B\'erczi-Szenes \cite{BercziSzenes12}. First, we explain a general framework.  
Let $V$ be a vector space. Suppose that $X$ is a compact algebraic manifold and $Y \subset X$ a subvariety. Let $E \to Y$ be a subbundle of the trivial bundle $X \times V \to M$ restricted to $Y$. Let $\pi_2: X \times V \to V$ be the projection, $i:E \subset X \times V$  the embedding,and $\phi=\pi_2\circ i$, as in the diagram

$$
\xymatrix{
E \ar@/^18pt/[rr]^{\phi} \ar@{^{(}->}[r]_-{i}\ar[d] & X \times V \ar[r]_-{\pi_2} \ar[d] & V\\
Y\ar@{^{(}->}[r] & X &\\
}
\eqno{(T)}
$$

Suppose that the torus $T$ acts on all spaces in the diagram above and that all maps are $T$-equivariant. 
Then, as a consequence from the Berline-Vergne-Atiyah-Bott Localization Formula \cite{AtiyahBott} applied to this setup, the following proposition holds. Here $e^T(-)$ is the Euler class, i.e., the top Chern class $c^T_{top}(-)$, 
and for a proper inclusion $\iota: A \hookrightarrow B$, we denote by $[A \subset B]$ the image of the fundamental class of $A$,  $\iota_*(1) \in H^*_T(B)$. 

\begin{prop} [B\'erczi-Szenes \cite{BercziSzenes12}] \label{localization}
Assume that the fixed point set $F(X)$ of the $T$-action on $X$ is finite. Then for the push-forward map $\phi_*:H_T^*(E) \to H_T^*(V)=H^*(BT)$, we have
$$\phi_*(1)=\sum_{f \in F(Y)}  \frac{[Y \subset X]|_f \cdot [E_f \subset V]}{e^T(T_fX)}.$$
Consequently, if $\phi$ is birational to its image, then the right-hand side of (9) is equal to $[\varphi(E)] \in H_T^*(V)=H^*(BT)$.
\end{prop}

If a fixed point $f$ is a smooth point of $Y$, then $e^T(T_fX) =e^T(T_fY) \cdot [Y \subset X]|_f$. 
Hence especially,  if $Y$ is non-singular, then 
$$\phi_*(1)=\sum_{f \in F(Y)} \frac{[E_f \subset V]}{e^T(T_fY)}.$$

Proposition \ref{localization} is applied to the following setting associated to a $\K$-singularity type $\eta$ of  finitely determined germs $(\C^m,0) \to (\C^{m+\kappa},0)$, according to Feh\'er-Rim\'anyi \cite{FR07b, FR12}. Namely, we consider    
$$V=J^k(m,m+\kappa)=J^k(m,1) \otimes \C^{m+\kappa}$$
and $X=\Gr^{\mu}$ is a suitable Grassmannian manifold including a certain subvariety $Y$, whose precise definition will be given later.  
This may be seen as an equivariant generalization of Damon's partial desingularization \cite{Damon}. 

Recall that the $\K$-equivalence classifies ideals $I_\eta$, the ideal generated by component functions of map-germs up to isomorphisms of the source space. Namely, a $\K$-type $\eta$ is determined by its local algebra $Q(\eta)=\Ost_{\C^m,0}/I_\eta$, so the dimension of the target space is somewhat auxiliary, i.e., $\kappa$ is a sort of parameter.  Since we are concerned with finitely determined germs $\eta$,  if $\kappa \ge 0$, the local algebra has finite dimension, and if $\kappa <0$,  it is enough to think of the local algebra modulo sufficiently high degree, i.e., $Q_k(\eta)=\Ost_{\C^m,0}/I_\eta+\mathfrak{m}_{\C^m,0}^{k+1}$, where $k$ is the determinacy degree of $\eta$. 

Let $\Hilb^{\mu+1}(\Ost_{\C^m, 0})$ denote the {\em punctual Hilbert scheme} which parameterizes ideals $I \subset \Ost_{\C^m, 0}$ supported at the origin of colength $\mu+1$ \cite{Iarrobino}. 
It is a projective variety, and especially, it is embedded in a Grassmannian $\Gr^{\mu}$ of $\mu$-codimensional subspace of the vector space  $J^k(m, 1)=\mathfrak{m}_{\C^m, 0}/\mathfrak{m}_{\C^m, 0}^{k+1}$. 
Given a $\K$-type $\eta$ with $\dim Q(\eta)=\mu+1$, we assign a stratum of $\Hilb^{\mu+1}(\Ost_{\C^m, 0})$ which  consists of all ideals being isomorphic to $I_\eta$. We denote the closure of the stratum by 
$$Y=Y_\eta \subset \Hilb^{\mu+1}(\Ost_{\C^m, 0}) \subset \Gr^{\mu}=X.$$ 

Note that $V, X, Y$ have actions of the maximal torus $T=(\C^*)^m$ of $GL_m$ in a natural way. 
Apply Proposition \ref{localization} to this setting. Then the information of the fixed points enables us to compute the $T$-equivariant Poincar\'e dual to the closure of the $\K$-orbit of type $\eta$ in $V=J^k(m,n)$, and hence, we get the Thom polynomial $tp(\eta)$ in $c_i=c_i(f)$ by the splitting principle. 

At the initial stage of this approach, B\'erczi-Szenes \cite{BercziSzenes12} considered the case of $\eta=A_{\mu}$ ($\kappa \ge 0$). 
In this case, $Y_{A_{\mu}}$ is realized as a particular component of the punctual Hilbert scheme, the {\em curvilinear Hilbert scheme}  $\CHilb^{\mu+1}(\Ost_{\C^m, 0})$, that consists of ideals isomorphic to $\C[x]/\langle x^{\mu+1} \rangle$, and it is a smooth quasi-projective variety. In fact, the authors dealt with  the {\em non-reductive group quotient} $J^k(1, m)\git \Aut(\C,0)$, via the action of the group of $k$-jets of re-parametrization germs of $\C$ at the origin, and compactify it by using the curvilinear Hilbert scheme. 
Their approach may be seen as a far generalization of Porteous-Gaffney's probe (test curve) method \cite{Porteous83, Gaffney83} from the equivariant geometry viewpoint. 
Their main result is the following iterated residue integral expression for the Thom polynomial of type $A_\mu$ where $\kappa \ge 0$ is regarded as a parameter. 

\begin{thm}[B\'erczi-Szenes \cite{BercziSzenes12}, B\'erczi \cite{Berczi17, Berczi20}] \label{thm_BS}
There exists a certain polynomial $Q_\mu$ as defined in \cite[\S 6]{BercziSzenes12} which depends only on $\mu$ and 
satisfies 
$$tp(A_\mu)=
 \underset{t_1=\infty}{\Res}\cdots \underset{t_\mu=\infty}{\Res}\, 
\frac{ \prod_{1\le i< j \le \mu} (t_j-t_i) \cdot {Q}_\mu(t_1, \cdots, t_\mu)}
{\prod_{k=1}^\mu \prod_{j=1}^{k-1}\prod_{i\le j}(t_k-t_i-t_j)}
\prod_{i=1}^\mu {\rm C}\hskip-2pt \left(\frac{1}{t_i}\right) t_i^\kappa\; dt_i 
$$
where ${\rm C}(t^{-1})=1+c_1t^{-1}+c_2t^{-2}+\cdots$ with the quotient Chern classes $c_i$. 
\end{thm}
Here, $ \underset{t_1=\infty}{\Res}\cdots \underset{t_\mu=\infty}{\Res} \Psi=(-1)^\mu \int_{|t_1|=R_1} \cdots \int_{|t_\mu|=R_\mu} \Psi \, dt_1\cdots dt_\mu$ where $0\ll R_1 \ll R_2 \ll \cdots \ll R_\mu$ for a rational function $\Psi(t_1, \cdots, t_\mu)$. 

\begin{exam}\upshape 
In B\'erczi-Szenes  \cite{BercziSzenes12}, $Q_\mu$ for $1\le \mu \le 6$ are explicitly given: 
$Q_1=Q_2=Q_3=1$ and $Q_\mu \, (\mu=4,5,6)$ are polynomials of degree $1, 3, 7$, respectively. 
For example, in case of $\mu=2$, we have 
$tp(A_2)$ as the integral of 
$$
\frac{t_2-t_1}{t_2-2t_1} \left(t_1^\kappa \sum_{i\ge 0} c_i t_1^{-i} \right)\left( t_2^\kappa \sum_{i\ge 0} c_i t_2^{-i}\right)\, dt_1dt_2
$$
and expanding it, we recover Ronga's formula (\ref{thom4}) written in \S \ref{fold_and_cusp}. 
For $\mu=3$,  we have $tp(A_3)$ by integrating 
$$
\frac{(t_2-t_1)(t_3-t_1)(t_3-t_2)}{(t_2-2t_1)(t_3-2t_1)(t_3-t_1-t_2)} \prod_{i=1}^3\left(t_i^\kappa \sum_{i \ge 0} c_i t_i^{-i} \right)\, dt_1dt_2dt_3
$$
which gives a simpler expression of a formula in B\'erczi-Feh\'er-Rim\'anyi \cite{BFR}. 
\end{exam}

Not only the local algebra of $A_n$-type but also the local algebra $Q(\eta)$ associated to any $\K$-type $\eta$ of finitely determined germs $(\C^m,0) \to (\C^{m+\kappa}, 0)$ with $\kappa\ge 0$, can be dealt with in the same framework. Combining it with the stabilization argument, the following strong consequence is obtained by Feh\'er and Rim\'anyi \cite{FR12}: 

\begin{thm}[Feh\'er-Rim\'anyi \cite{FR12}] \label{thm_Ts} 
{\bf (Thom series)} 
Let $Q$ be a local algebra of finite colength $\mu+1$, then there is a unique formal series 
$${\rm Ts}_Q=\sum_{i_1, \cdots, i_\mu} \psi_{i_1,\cdots, i_\mu} c_{\kappa+1+i_1} \cdots c_{\kappa+1+i_\mu}$$
such that for any $\K$-type $\eta$ of map-germ with $\kappa \ge 0$ whose local algebra $Q(\eta)$ is isomorphic to $Q$, 
the corresponding Thom polynomial $tp(\eta)$ is obtained as ${\rm Ts}_Q$ evaluated by the quotient Chern classes $c_{\kappa+1+i}=c_{\kappa+1+i}(f)$ with conventions that $c_0=1$ and $c_j=0$ for $j<0$. Here the integer coefficients $\psi_{i_1,\cdots, i_\mu}$ are independent from $\kappa$. 
\end{thm}

In fact, ${\rm Ts}_Q$ formally consists of infinitely many terms indexed by $\mu$-tuples of integers $(i_1,\cdots, i_\mu)$ whose sum is some constant determined by $Q$, but fixing $\kappa$ and applying the conventions on $c_j$, it is read off as a polynomial ($=tp(\eta)$). This is the most comprehensive existence theorem of universal polynomials for mono-singularity types -- ${\rm Ts}_Q$ is called the {\em Thom series for the local algebra $Q$}. 
For a nice introduction to Thom series with several examples, see Rim\'anyi's survey article \cite{Rimanyi24}. 

The rule of the appearance of the coefficients $\psi_{i_1,\cdots, i_\mu}$ should be detected. 
It requires to replace the Grassmannian $\Gr^n$ by a finer space. In \cite{Kaz17} Kazarian has introduced the concept of {\em non-associative Hilbert scheme} to be the total space of some (partial) flag bundle which becomes a smooth ambient space into which $Y_\eta$ is embedded as a closed subvariety. Let $Q=Q(\eta)$ be the local algebra of a $\K$-type $\eta$ of map-germs and let  $\dim Q=\mu+1$. In  order to construct this ambient space, a key ingredient is a filtration on the nilpotent subalgebra $N:=\mathfrak{m}_{\C^m,0}Q$, 
$$N=N_1 \supset N_2 \supset \cdots \supset N_{r+1}=0 \qquad d_i=\dim (N_i/N_{i+1})$$
satisfying $N_i\cdot N_j \subset N_{i+j}$. Of course, $\sum d_i=\mu$. 
The degree filtration with $N_k:=\mathfrak{m}^k_{\C^m,0}Q$ is a typical one, but there may be some other choices. 
We call the sequence $\underline{d}=(d_1, \cdots, d_r)$ the {\em dimension vector} associated to the filtration on $N$. 
Put $s_k=d_1+\cdots+d_k$, and for each $i$ with $s_{k-1} \le i \le s_k$, put $e(i):=s_k-i$ and $w(i)=k$, and set 
$$K_{\underline{d}}(t_1, \cdots, t_\mu)
:= \frac{\prod_{i=1}^\mu t_i^{e(i)} \cdot \prod_{1\le i< j \le \mu} (t_j-t_i) }
{\prod_{k=1}^d \prod_{j=1}^{k-1}\prod_{i\le j, w(i)+w(j) \le w(k)}
(t_k-t_i-t_j)}.$$

\begin{thm}[Kazarian \cite{Kaz17}] \label{Kaz_Ts} 
The coefficients $\psi_{i_1,\cdots, i_\mu}$ arising in ${\rm Ts}_Q$ are determined by the following rational generating function
$$\sum_{i_1, \cdots, i_\mu} \psi_{i_1,\cdots, i_\mu} t_1^{i_1}\cdots t_\mu^{i_\mu} 
= K_{\underline{d}}(t_1, \cdots, t_\mu) P_Q(t_1, \cdots, t_\mu)$$
where $P_Q$ is a polynomial defined as in \cite[\S 4]{Kaz17}.  
In other words, the Thom series ${\rm Ts}_Q$ has the following residual integral expression: 
$${\rm Ts}_Q=\underset{t_1=\infty}{\Res}\cdots \underset{t_\mu=\infty}{\Res} \; 
K_{\underline{d}}(t) \cdot P_Q(t) \cdot \prod_{i=1}^\mu {\rm C}(t_i^{-1})\, t_i^\kappa \; dt_i.$$
\end{thm}

Obviously, in the case of $\eta=A_\mu$ with the degree filtration, this residual integral is just Theorem \ref{thm_BS}. 
Besides from $A_\mu$-singularity types, for other types, e.g., $\Sigma^{i, j}$, $I_{a, b}$, $I \hskip-2pt I \hskip-2pt I_{a,b}$ etc., the polynomial $P_{Q(\eta)}$, and hence ${\rm Ts}_{Q(\eta)}$, are explicitly computed, see \cite{Kaz17, FR12}. 

Note that this decomposition (the right hand side) of the generating function of $\psi_{i_1,\cdots, i_\mu}$ is not unique, especially depends on the choice of filtration; $K_{\underline{\eta}}$ is the factor coming from  the Gysin map of the projection of the flag bundle (cf. \cite{Damon73}) and $P_Q$ corresponds to the fundamental class of a certain Borel orbit in the ambient space. 

Looking back the flow of old and new researches, we observe that the present work of Thom series surely lies on the natural extension of work of Porteous, Ronga, Damon, Gaffney etc. That has not yet been completed, rather say, still only halfway there. In this theory, it is most critical to detect the polynomial $P_Q$, which comes from the Borel geometry associated to the local algebra, and it essentially relates to tautological integral over punctual Hilbert schemes.  Many tasks remain for future work.

\section{Higher Thom polyomials}\label{higher_tp}
The Thom polynomial $tp(\eta)(c(f))$ is the most basic invariant for the appearance of $\eta$-type singular points for a given map $f: M \to N$, while if the locus $\overline{\eta(f)}$ has the expected dimension which is positive,  the locus may contain higher codimensional informative cycles other than the fundamental cycle. 
For instance, if $f: M \to N$ is a Morin map  (a stable map which has only local singularities of type $A_k$), the locus $S_k:=\overline{A_k(f)}$ is a submanifold of $M$,  and hence, the pushforward of the Chern class $\iota_*c(TS_k) \in H^*(M)$ is defined ($\iota$ is the inclusion of the locus to $M$). These classes appeared in  \cite{HaroldLevine, Ando, Nakai}, and especially, Ando called them `higher Thom polynomials for Morin singularities'. 
However, the locus for non-Morin maps is singular in general, then this definition does not make sense at all. We need some more tools to build a theory.

\subsection{Chern class for a singular variety} 
Let $X$ be an irreducible algebraic variety of dimension $n$. 
If $X$ is non-singular, it possesses the {\em total} Chern class 
$$c(TX)=1+c_1(TX)+\cdots + c_n(TX) \; \in H^*(X).$$
Switching to homology, note that $c_n(TX) \frown [X] = \chi(X)$ and $1\frown [X]=[X]$. 
If $X$ is singular, the tangent bundle no longer exists, so $c(TX)$ does not make sense at all. 
Nevertheless, we may find some alternative to $c(TX)$ in homology -- for example, 
the {\it Chern-Schwartz-MacPherson class} (CSM class) is a typical one, see Brasselet \cite{Brasselet22}. 
It is a total homology class of $X$ such that 
the degree is topological Euler characteristic of $X$ and the top term is the fundamental class
$$c^{\SM}(X)=\chi(X) \cdot [pt] + \cdots + [X] \; \in H_*(X)$$
together with a nice functorial property (let us note again that if $X$ is non-compact, $H_*$ stands for 
Borel-Moore homology).  
Let $\F(X)$ be the abelian group of constructible functions on $X$, and 
for a proper morphism $X\to Y$, the {\it pushforward} $f_*: \F(X)\to \F(Y)$ is defined. 
It holds that $(f\circ g)_*=f_*\circ g_*$. 
R. MacPherson introduced the following natural transformation between two covariant functors $\F$ and $H_*$ 
using the Chern-Mather class and local Euler obstruction, as a solution to the so-called Grothendieck-Deligne conjecture. 
Also it should be mentioned that the same notion had been found by M. Schwartz in a different context from obstruction theory \cite{Brasselet22}.

\begin{thm}\label{C_*} \cite{Mac}
There is a unique natural transformation 
$$C_* : \F(X) \longrightarrow H_*(X)$$ 
between these functors 
so that $C_*(\jeden_X)=c(TX)\frown [X]$  if $X$ is non-singular. 
\end{thm}
Namely, it holds that  
\begin{itemize}
\item $C_*(\alpha+\beta)=C_*(\alpha)+C_*(\beta)$ (additive homomorphism) 
\item 
$C_*f_*(\alpha)=f_*C_*(\alpha)$ for proper morphisms $f: X \to Y$. 
\end{itemize}
In particular,  
if $pt: X\to pt$ is proper, we have 
$pt_*C_*(\alpha)=C_*pt_*(\alpha)=\int_X \alpha$.  
The {\it Chern-Schwartz-MacPherson class} of $X$ 
is defined by $c^{\SM}(X):=C_*(\jeden_X)$.


There is also the `Segre-version' of CSM classes, see  Aluffi \cite{Al22, Al05}. 
Let $\iota: X \hookrightarrow M$ be a closed embedding into non-singular $M$. 
Like as the relation between Segre covariance class and Fulton's Chern class, 
we define the {\it Segre-Schwartz-MacPherson class} by 
$$s^{\SM} (X,M):= c(\iota^*TM)^{-1}\frown  c^{\SM}(X) \;\; \in H_*(X).$$
We often regard the class $s^{\SM} (X,M)$ as in $H^*(M)$ via the pushforward $\iota_*$ and the Poincar\'e duality, unless specifically mentioned. 
Also we set for $\alpha \in \F(M)$ 
$$s^{\SM} (\alpha,M) := c(TM)^{-1}\cdot  C_*(\alpha)\in H^*(M).$$  
The Segre-SM class has an expected nice property for transverse pullback: 

\begin{prop}\label{segre_pullback} 
Let  $f:M \to N$ be a map between complex manifolds, and 
let $Y$ be a closed singular subvariety of $N$.  
Assume that $f$ is transverse to (a Whitney stratification of) $Y$. 
Then it holds that 
$$f^*s^{\SM}(Y,N)=s^{\SM}(f^{-1}(Y),M) \;\;\;\in \;\;H^*(M).$$
\end{prop} 

\begin{rem}\label{sign_convention}\upshape 
We follow the sign convention of the Segre class as in Fulton's book \cite{Fulton}. 
Notice that if $X$ is a closed submanifold of $M$ with the normal bundle 
$\nu = \iota^*TM - TX$,  
then the Segre-SM class is nothing but 
the {\it inverse normal Chern class} for $X \hookrightarrow M$: 
$$s^{\SM} (\jeden_X,M)=\iota_* c(-\nu) \; \; \in H^*(M).$$ 
\end{rem}

Based on equivariant intersection theory \cite{Totaro, EdidinGraham}, an equivariant theory of CSM class is available. 
Let $\F^G_{inv}(X)$ be the subgroup generated by $G$-invariant constructible functions in $\F(X)$. 
Using $G$-equivariant Chern-Mather class and local Euler obstruction, one gets the following theorem\footnote{In \cite[p.122]{Ohmoto06}, another construction of $C_*^G$ using inductive limit is also given. Note that if $G$ is not a torus or a finite group, a tiny modification is needed to multiply $T_{U, *}$ by a correction term to guarantee the normalization condition. This does not affect the correctness of Theorem \ref{equiv_CSM}, though.}.

\begin{thm}  \label{equiv_CSM}
\cite{Ohmoto06, Ohmoto12a, Ohmoto16, Weber12} 
For $G$-varieties and proper $G$-morphisms, 
there is a unique natural transformation 
$$C_*^G: \F^G_{inv}(X) \to H^G_*(X)$$ 
so that $C_*^G(\jeden_X)=c^G(TX) \frown [X]_G$ if $X$ is non-singular. 
\end{thm}

Each dimensional component of the equivariant CSM class 
has its support on an invariant algebraic cycle in $X$. 
In particular, the lowest and highest terms are as follows: 
If $X$ is of equidimension $n$, the top term is the $G$-equivariant fundamental class:  
$C_n^G(\jeden_X) = [X]_G$. 
If $X$ is compact, the degree $pt_*^G C_0^G(\alpha)$ is equal to the weighted Euler characteristics $\int_X \alpha$, especially  
$pt_*^G C_0^G(\jeden_X)=\chi(X)$. 

Also the Segre-SM classes admit an equivariant verision: for $G$-embedding $X \hookrightarrow V$ where $V$ is a $G$-manifold, we define 
$$
s_G^{\SM}(X, V):=c^G(TV)^{-1} \frown C^G_*(\jeden_X) \; \in H_G^*(X)
$$
where $c^G(TV)$ is the $G$-equivariant Chern class of $TV$. 
Abusing the notation, we also mean its pushforward into $M$ and take the equivariant Poincar\'e dual: 
$s_G^{\SM}(X, V) \in H_G^*(V)$ unless any confusion occurs. 
It immediately follows from Proposition \ref{segre_pullback} that the equvariant SSM class admits the naturality for the transversal pullback via $G$-morphism. 

In particular, if $V$ is a vector space with $G$-action and $W \subset V$ is a $G$-subvariety, we define 
$$tp^{\SM}(W) \, (\, = tp^{\SM}(\jeden_W)\,)\,:=s_G^{\SM}(W, V) \;\;\;  \in H_G^*(V)\simeq H^*(BG).$$
Its leading term is just the Thom polynomial $tp(W) \in H^*(BG)$. 
We call $tp^{\SM}(W)$ the {\it higher Thom polynomial of $W \subset V$} with respect to the SSM class. 
This is a universal SSM class in the following sense. 

\begin{thm} \label{univ_ssm}
\cite{Ohmoto06, Ohmoto12a, Ohmoto16} 
Given a vector bundle $\pi: E \to M$ over a complex manifold $M$ 
with fiber $V$ and structure group $G$,  
let  $W(E) \to M$ be the associated cone bundle with the fiber $W$ and group $G$. 
Then, for any section $s: M \to E$ which is transverse to $W(E)$, 
it holds that 
$$\Dual\, s^{\SM}(W(s), M)=\rho^* tp^{\SM}(W) \in H^*(M)$$
where $W(s):=s^{-1}(W(E))$ and $\rho$ is the classifying map for  $E \to M$. 
The topological Euler characteristic of $W(s)$ is given by 
$$\chi(W(s))=\int_M c(TM)\cdot \rho^*tp^{SM}(W).$$ 
\end{thm}

\begin{rem}\upshape \label{PP}
{\bf (CSM class of degeneracy loci)} 
A prototype of Theorem \ref{univ_ssm} is seen in Parusinski-Pragacz \cite{PP}, see also in Fulton-Pragacz \cite{FP98}.  They actually considered the CSM class $c^{\SM}(\overline{\Sigma^k(\varphi)})$ of degeneracy loci of vector bundle maps $\varphi: E \to F$ over $M$ by using the desingularization method, as the CSM class generalization of the Thom-Porteous formula. 
More precisely, for all $s \ge k$, they collect Porteous' desingularizations $\pi_s: Z_s \to \overline{\Sigma^s(\varphi)}$ mentioned in \S \ref{desing} and compute $\pi_{s*}c^{\SM}(Z_s) \in H^*(M)$ in terms of Chern classes of $E$ and $F$. It then follows from the additivity (deletion-inclusion principle) of CSM classes to see that  
$$c^{\SM}(\overline{\Sigma^k(\varphi)})=\pi_{k*}c^{\SM}(Z_k)+\sum_{s>k} a_s \pi_{s*}c^{\SM}(Z_s)$$ 
for some $a_s \in \Z$, and those coefficients $a_s$ are computed using Euler characteristics of several Grassmannians (singular fiber of $\pi_s$'s). That yields Parusinski-Pragacz's formula. 
\end{rem}

\begin{rem}\upshape \label{GRR}
{\bf (Grothendieck-Riemann-Roch)}
Note that $C_*$ can be thought of as a {\em singular} variant of Grothendieck-Riemann-Roch type theorem (GRR) in Grothendieck's original sense. 
The most general natural transformation in a certain sense is the {\em Hirzebruch class transformation} given by Brasselet-Sch\"urmann-Yokura \cite{BSY}, which unites simultaneously  

\begin{itemize}
\item[-] MacPherson's Chern class transformation, 
\item[-] Baum-Fulton-MacPherson's Todd class transformation (singular GRR),  
\item[-] Cappell-Shaneson's $L$-class transformation
\end{itemize}
(see surveys in Handbook \cite{Brasselet22, Yokura} for more information). 
$G$-equivaraint version is also available through the algebraic Borel construction -- especially, using localization technique, the torus-equivariant Hirzebruch class transformation has been studied by Weber \cite{Weber13, Weber16, Weber17}, see also \S  \ref{modern_Schubert} below. 
\end{rem}

\subsection{Universal Segre-Schwartz-MacPherson classes}
Let us apply Theorem \ref{univ_ssm} to the classification of linear maps (\S \ref{thom_porteous}) and the $\K$-classification of singularities of maps (\S \ref{fold_and_cusp}). 

\subsubsection{Thom-Porteous formula, upgraded} \label{modern_Schubert} 
This research direction has been upgraded into modern Schubert calculus and geometric representation theory. 

Based on the formulation in \cite{Ohmoto06, Weber12}, 
Feh\'er-Rim\'anyi \cite{FR18} has developed a computational and representation-theoretic aspect on equivariant CSM/SSM classes for orbit-closures of $G$-action, e.g., generalized Thom-Porteous formula for $\Sigma^k$-types (cf. Remark \ref{PP}) -- in particular, they provide the interpolation characterization. 
Independently from this work, a different approach to (equivariant) CSM/SSM class of Schubert varieties was initiated by  Aluffi-Mihalcea \cite{AM09, AM16},  and a more comprehensive study was achieved by Aluffi-Mihalcea-Sch\"urman-Su \cite{AMSS23, AMSS22}. 

In particular, preceding to \cite{FR18}, the torus-equivariant CSM class has been featured by Rim\'anyi-Varchenko \cite{RV15} in the context of cotangent Schubert calculus, where they use it to identify {\em cohomological stable envelope} in the sense of A.~Okounkov. This is a striking result, indeed. It is soon lifted up to {\em $K$-theoretic stable envelope} -- the counterpart is the {\em equivariant motivic Chern classes}, i.e., equivariant version of Brasselet-Sch\"urmann-Yokura's motivic Chern classes \cite{BSY}. In recent years, this direction has been intensively studied by Rim\'anyi \cite{Rimanyi20, Rimanyi21}, Feh\'er \cite{Fe21}, Feh\'er-Rim\'anyi-Weber \cite{FRW17, FRW21}, and also Aluffi-Mihalcea-Sch\"urman-Su \cite{AMSS22, AMSS23}, and moreover, for an elliptic cohmology version of stable envelopes, by Rim\'anyi-Weber \cite{RW20, RW22}. For the relation with quantum integrable systems and mirror symmetry, see \cite{RSVZ19, RSVZ22, RS22, RS22} and related references. 
The theory is still developing rapidly (several new papers have been produced continuously, e.g.,  \cite{KRW20, KW22, KW23, PR22}), and therefore, it is well beyond the scope of  this survey paper.  For the starting point of this development, the readers should be referred to,  e.g., Rim\'anyi \cite{Rimanyi21}, Feh\'er-Rim\'anyi-Weber \cite{FRW21}, Aluffi-Mihalcea-Sch\"urman-Su \cite{AMSS23} and references therein.

\subsubsection{Fold and Cusp, upgraded} \label{fold_and_cusps_upgraded} 
Let $\eta$ be a singularity type in $\K$-classification of map-germs with codimension $\kappa$. 
Now we apply Theorem \ref{univ_ssm} to the situation with $G=J^k\K_{m,m+\kappa}$ acting on the jet space $V=J^k(m,m+\kappa)$ and  the closure $\overline{\eta} \subset V$. Then by entirely the same proof of Theorem \ref{tp}, we see that there is a unique universal power series 
$$tp^{\SM}(\overline{\eta})=tp(\eta)+h.o.t  \;\;\; \in \Z[[c_1,c_2, \cdots ]]$$
so that for any generic map 
$f : M \to N$ of  codimension $\kappa$ 
it holds that 
$$
\Dual s^{\SM}(\overline{\eta(f)}, M) =  tp^{\SM}(\overline{\eta})(c(f)) \;\; \in H^*(M)
$$
after substitution by the quotient Chern classes $c_i=c_i(f)$  \cite{Ohmoto06, Ohmoto16}.

\begin{exam}\label{exam_tpsm}
\upshape 
The interpolation method of Rim\'anyi is also very effective to compute $tp^{\SM}(\overline{\eta})$ as well as the leading term (Thom polynomial)  $tp(\eta)$. 
For stable singularity type $\eta$ of maps in case of $\kappa=0$, we compute 
the truncated polynomials of $tp^{\SM}(\overline{\eta})$ up to degree $4$ as follows: 
\begin{eqnarray*}
tp^{\SM}(\overline{A_1})&\equiv& 
c_1-c_1^2+c_1^3-c_1^4+c_2^2-c_1c_3   \\
tp^{\SM}(\overline{A_2}) &\equiv& 
c_1^2+c_2 - (2c_1^3+3c_1c_2+c_3)+3c_1^4+6c_1^2c_2+4c_2^2+c_4-6c_2^2+6c_1c_3   \\
tp^{\SM}(\overline{A_3}) &\equiv& 
c_1^3+3c_1c_2+2c_3-(3c_1^4+12c_1^2c_2+15c_2^2+6c_4)+14c_2^2-14c_1c_3   \\
tp^{\SM}(\overline{A_4}) &\equiv& 
c_1^4+6c_1^2c_2+2c_2^2+9c_1c_3+6c_4  \\
tp^{\SM}(\overline{I_{2,2}}) &\equiv& c_2^2-c_1c_3  
\end{eqnarray*}
Furthermore, in cases of $0\le \kappa \le 4$,  $tp^{\SM}$ for $\K$-types are explicitly computed up to certain codimension, see Rim\'anyi's 'Thom polynomial Portal' \cite{RimanyiPortal}. 
In fact,  based on the previous work with Feh\'er \cite{FR18} for SSM classes of degeneracy loci, Rim\'anyi \cite{Rimanyi25} has  developed the interpolation for computing $tp^{\SM}$ in connection with geometric representation theory of Maulik-Okounkov's stable envelopes (cf. \S \ref{modern_Schubert}). 
Those higher Thom polynomials are useful, e.g., in classical enumerative geometry: 
if some $\eta$-singular locus is a curve, then we can compute the degree and the genus of its normalization by using the first two terms of $tp^{\SM}(\overline{\eta})$, and 
those information is needed to recover classical enumerative formulas \cite{Sasajima17, Nekarda}. 
\end{exam}

Like as this example, we could do compute $tp^{\SM}(\overline{\eta})$ up to certain degree term, but it is not realistic to compute all terms, because, as a theoretical reason, there appear infinitely many strata adjacent to $\eta$ under the stabilization process ($m \to \infty$), and as a practical reason, if one uses the interpolation method to compute the SSM class,  the same difficulty arises as when computing $tp$, namely  moduli of $\K$-orbits appear at certain codimension. 

Nevertheless, we may have a chance to find a closed expression of  $tp^{\SM}$ for some nice invariants. Suppose that we are given a proper stable map $f: M \to N$ with $d=\dim M - \dim N (=-\kappa) \ge 0$. 
Note that the germ $f:(M, x) \to (N, f(x))$ at every critical point $x$ defines an $d$-dimensional isolated complete intersection singularity (ICIS). 
Then we can define a constructible function over $M$, called 
the {\it Milnor number constructible function} associated to $f$: 
$$\mu(f): M \to \Z, \quad x \mapsto \mbox{Milnor number of $f$ at $x$}.$$
In fact, $\mu: J^k(m,n)-W_k \to \Z$, where $W_k$ is some $\K$-invariant closed algebraic set, is defined as a $\K$-invariant constructible function such that $\mu(f)$ is obtained by the pulling back $\mu$ via the jet extension and the codimension of $W_k$ goes to $\infty$ ($k \to \infty$) \cite[\S 7]{Gibson}. Thus a universal series $tp^{\SM}(\mu)$ exists. 
On the other hand, for a proper  stable map $f: M \to N$ with connected $N$ and  generic fibre $F$,  
it holds that 
$$f_*(\, \jeden_M+(-1)^{d}\mu(f) \, )= \chi(F) \cdot \jeden_N$$
(the Yomdin-Nakai formula \cite{Yomdin, Nakai}). 
Apply $C_*$ to the both sides of this equality and multiply by $c(TN)^{-1}$, and then, after a small calculation, 
universal relation among Landweber-Novikov classes $s_I$ suggests the following conjectural formula \cite{Ohmoto12b, Nakai}: 
\begin{equation}\label{Milnor1}
s^{\SM}(\mu(f), M)= tp^{\SM}(\mu) (c(f))= (-1)^{d}\, \left(\, \sum_{i=0}^\infty \sum_{j=0}^d c_i(f)\bar{c}_j(f)-1\, \right) 
\end{equation}
in $H^*(M)$. That is, 
\begin{eqnarray*}
tp^{\SM}(\mu) 
&=&(-1)^d ((1+c_1+c_2+\cdots)(1+\bar{c}_1+\cdots + \bar{c}_d)-1) \\
&=&(-1)^{d+1}(1+c_1+c_2+\cdots)(\bar{c}_{d+1}+\cdots).
\end{eqnarray*}
There is an attempt to give a direct proof using conormal geometry by X. Liao \cite{Liao23}. 

\begin{exam}\upshape 
In particular, in case of $\dim M=d+1$ and $\dim N=1$, the above (\ref{Milnor1}) recovers the following classical formula  (cf. Fulton \cite{Fulton}):  
if $f$ has finitely many isolated critical points, then 
the sum of the Milnor number is expressed by
$$
 \int_M \mu(f)   =  (-1)^{d +1} c_{d +1}(TM-f^*TN)\frown [M]. 
$$
As another example, in case of $d\, (=-\kappa)=0$, 
we see 
$\mu=
\jeden_{\overline{A_1}}+\jeden_{\overline{A_2}}+\jeden_{\overline{A_3}}+\jeden_{\overline{A_4}}+\cdots 
\in \F^{\K}(J(n,n))
$ 
modulo constructible functions with support of codimension greater than $4$. 
Summing up $tp^{\SM}(\overline{\eta})$ in Example \ref{exam_tpsm},  we observe, up to degree four, 
$tp^{\SM}(\mu)\equiv tp^{\SM}(\overline{A_1}) + \cdots + tp^{\SM}(\overline{A_4}) 
\equiv c_1+c_2+c_3+c_4$. This is consistent with $tp^{\SM}(\mu)=-(1+c_1+\cdots)(\bar{c}_1+\cdots)=-c(\bar{c}-1)=c_1+c_2+c_3+c_4+\cdots$ 
(cf.  \cite[p.223]{Ohmoto12a}). 
\end{exam}

\begin{exam}\upshape \label{icis}
{\bf (Milnor number of quasi-homogeneous ICIS)} 
Suppose that we are given a quasi-homogeneous ICIS 
$$\eta: (\C^m,0) \to (\C^n,0) \qquad (d=m-n\ge 0)$$ 
with weights $\bw=(w_1, \cdots , w_m)$ and degrees $\bd=(d_1, \cdots, d_n)$. 
Localizing the above conjectural formula (\ref{Milnor1}) via the action of $T=\C-\{0\}$, we get a new interpretation or better understanding of a well-known Milnor number formula of Greuel-Hamm and Damon \cite{Damon96}. 
In the same way as in \S \ref{fold_and_cusp2} we associate the universal map $f_\eta: E_0 \to E_1$, then   
$c(E_0)=\prod_{i=1}^m (1+w_i t)$, $c(E_1)=\prod_{j=1}^n (1+d_j t)$ 
and $c(E_0-E_1)=c(E_0)/c(E_1)$. 
Applying $tp^{\SM}(\mu)$ to $f_\eta$, we obtain an expression of the Milnor number of $\eta$: 
\begin{equation} \label{Milnor2}
\mu_\eta = (-1)^{d}\left(\frac{c_n(E_1)}{c_{m}(E_0)}c_{d}(E_0-E_1)-1\right). 
\end{equation}
This is actually the same as the formula of Greuel-Hamm obtained in an entirely different way: 
$$\mu_\eta=(-1)^{d+1}+\frac{|\bd |}{| \bw |}\sum_{j=0}^{d} 
(-1)^{d-j}\sigma_{d-j}(\bw)s_j(\bd)$$
where $|\bw|=\prod_{i=1}^m w_i$, $|\bd|=\prod_{i=1}^n d_i$ and 
$\sigma_k(\bw)$ and $s_k(\bd)$ are respectively 
the $k$-th elementary symmetric function in $w_1, \cdots, w_m$ and 
the $k$-th complete symetric function in $d_1, \cdots, d_n$, i.e., Chern (resp. Segre) classes. 
\end{exam}

\section{Other classifications of map-germs}\label{variants}
We give a quick review on Thom polynomial theory for other classifications of map-germs, especially, 
isolated singularities of functions and hypersurfaces, and $\A$-finitely determined map-germs.

\subsection{Lagrange and Legendre Thom polynomials}\label{legendre}
There is an important class of mappings, called {\em Lagrange and Legendre maps}, whose singularities correspond to singularities of functions and hypersurfaces, respectively. 

Let  $Q$ be a complex symplectic manifold and $\pi: Q \to N$ a Lagrange fibration (i.e., every fiber is Lagragian), and 
$\iota: L \to Q$ a Lagrange immersion.  
A typical example is $Q=T^*N \to N$, the cotangent bundle equipped with the canonical symplectic structure. 
The composed map $\pi\circ \iota: L \to N$ is called a Lagrange map; locally such a map is expressed by the gradient maps $(\C^n,0) \to (\C^n, 0)$, $p \mapsto \frac{\rd g}{\rd p}$, of some generating function $g=g(p)$. 
There is the counterpart in contact geometry; the composition of a Legendre immersion and the projection of Legendre fibration is a Legendre map. 
Their singularities are called  {\em Lagrange and Legendre singularities}, respectively, and have thoroughly been explored by Arnold, Zakalyukin and others, see \cite{AGV}. The corresponding Thom polynomial theory was initiated by Vassiliev in real $C^\infty$-category \cite{Vassiliev, AGLV}, in which the simplest Thom polynomial for fold singularity corresponds to the so-called the (first) {\em Maslov class}. The theory is also formulated in complex analytic/algebraic category by Kazarian \cite{Kaz98, Kaz00, Kaz03a}, that we will briefly review below. Also see Mikosz-Pragacz-Weber  \cite{MPW09, MPW11}) especially for the positivity question. 
Some other particular cases had been treated in e.g.,  \cite{Hayden, MSV89}.

A hypersurface-singularity of dimension $d$ is a $\K$-equivalence class of function-germs $f: (\C^{d+1},0)\to (\C,0)$. 
It is an isolated singularity if and only if $f$ is $\K$-finitely determined. 
When ignoring the linear action on the target $\C$, it is the $\mathcal{R}$-classification of function-singularities. The list of typical singularity types starts from the so-called {\em ADE classification}. Although it might be thought of as a special part (with corank one) of the classification of isolated complete intersection singularities (ICIS),  it has own particular nature -- there is a different stabilization.  Namely, $f: (\C^{d+1},0)\to (\C,0)$ and $g: (\C^{d'+1},0)\to (\C,0)$ ($d\le d'$) are {\em stablly $\mathcal{G}$-equivalent} if $f(x)+x_{d+2}^2+\cdots + x_{d'+1}^2$ is $\mathcal{G}$-equivalent to $g$ ($\mathcal{G}=\mathcal{R}$ or $\K$). 

Below we focus mainly on the Legendrian case ($\mathcal{G}=\K$). 
A natural way to study global properties of isolated hypersurface singularities is described as follows. 
That arises in a family of divisors. 
Suppose that we are given a locally trivial complex fibration $p: W \to N$ where $W$ and $N$ are complex manifolds of dimension $m+1$ and $n$, respectively, and an embedding $\iota: M \hookrightarrow W$ of a non-singular hypersurface. Put 
$$\varphi=p\circ \iota: M \to N$$
($\dim M=m$, $\dim N=n$) and  $L$ to be the critical locus of $\varphi$, i.e., the locus of points $x \in M$ at which $T_xM$ contains the relative tangent space $\ker dp_x$. Note that the map $\varphi$ has singularities of corank at most one: 
$$L=\Sigma^{m-n+1}(\varphi), \quad \Sigma^{m-n+i}(\varphi)=\emptyset  \;\;\; (i>1).$$ 
In other words, the germ of $\varphi$ at each critical point is an unfolding of a function-germ. 
We now assume that $M$ is sufficiently generic so that the $1$-jet exptension $d\varphi$ is transverse to the first Thom-Boardman stratum $\Sigma^{m-n+1}$, then $L$ becomes a manifold of dimension $n-1$. In particular, the restricted map 
$$\varphi|_L: L \to N$$ 
is regarded as a Legendre map. Indeed, for $\dim \coker d\varphi_x=1$ at every $x \in L$,  we assign a covector defining the image hyperplane in $T_{\varphi(x)}N$ (up to non-zero constant multiples), then it turns out that the map $\varphi|_L$ is factored via a Legendre immersion:  
$\varphi|_L: L \to PT^*N \to N$.

Now we put $V=\ker d\varphi$ and $I=\coker d\varphi$; those are vector bundles over $L$ of rank $m-n+1$ and of rank one, respectively. 
The map $\varphi$ yields the germ of a non-linear map $V \to I$ at the zero section over $L$. Let $f$ denote the map-germ: 
$$
\xymatrix{
V \ar[rd] \ar[rr]^f& & I \ar[ld] \\
&L&
}
$$
The jet $j^kf(x)$ defines a section of 
$$J:=\bigoplus_{k\ge 2} {\rm Sym}^k V^* \otimes I.$$ 
This section $L \to J$ is actually determined by the $k$-th intrinsic derivatives of $\varphi$ ($k \ge 2$),  
and is most essential in our argument.

Let $\eta$ be a stable $\K$-type of isolated hypersurface singularities. Put $\eta(M) \subset L$ to be the locus consisting of the points $x \in L$ at which $M$ is tangent to the fiber with contact type $\eta$, in other words, the $\eta$-singularity locus of the map $\varphi: M \to N$. 
Hence, if $M$ is appropriately generic,  $[\eta(M)]$ is expressed  in $H^*(M)$ by the Thom polynomial $tp(\eta)$ in the quotient Chern classes $\bar{c}_i=c_i(TM-\varphi^*TN)$. However, the geometric setting makes a strong restriction, and especially,  $\Dual [\eta(M)]$ should be localized into $H^*(L)$. 
We are interested in representing the class using $V, I$.  
More precisely, we start with the above diagram and put
$$u= c_1(I), \quad a_i=c_i(V^*\otimes I-V),$$ 
where these classes satisfy relations
$$(1+a_1+a_2+\cdots ) \left(1-\frac{a_1}{1+u}+\frac{a_2}{(1+u)^2}-\cdots \right) =1.$$
Note that when writing $U=V^*\otimes I -V$ formally,  this relation follows from the identity 
$$U^*\otimes I+U= 0.$$ 
That allows us to expand the squares of the classes $a_i$, and hence any polynomial in $u, a_1, a_2, \cdots$ can be expressed as a linear combination of monomials $u^{i_0}a_1^{i_1}a_2^{i_2}\cdots $ with $i_0 \ge 0$ and $i_k=0$ or $1$ for $k >0$. 
We define the ring $\mathcal{L}$ of  {\em universal Legendre characteristic classes} to be the quotient ring of the polynomial ring $\Z[u, a_1, a_2, \cdots]$ with $\deg u=1$ and $\deg a_i=i$ by the ideal obtained by the above relations. 

\begin{thm}\label{thm_tp_L} 
{\bf (Legendre Thom polynomials)} \cite{Kaz03a} 
Let $\eta$ be a stable $\K$-type of isolated hypersurface singularities. 
There is a unique polynomial $P_\eta(u, a_i)$ in $\mathcal{L}$ such that 
for any generic map $f: V \to I$ over $L$ as above, i.e., any generic section $jf: L \to J$,  it holds that 
$$\Dual \,[\overline{\eta(f)}] = P_\eta(u, a_i) \; \in H^*(L)$$
evaluated by $u= c_1(I)$ and $a_i=c_i(V^*\otimes I-V)$. 
\end{thm}

In \cite{Kaz03a},  $P_\eta$ are explicitly computed for $\eta$ of codimension up to $6$, e.g., 
a few sample of them are written below (see also Remark \ref{lagrange}). 

\begin{align}
P_{A_2} &= \bm{a_1}\notag \\
P_{A_3} &= \bm{3a_2}+ ua_1  \notag \\
P_{A_4} &= \bm{3a_1a_2+6a_3}+4ua_2+u^2a_1  \notag \\
P_{D_4} &= \bm{a_1a_2-2a_3}-ua_2  \notag \\
P_{A_5} &= \bm{27a_1a_3+6a_4}+u(16a_1a_2-12a_3)-4u^2a_2+u^3a_1  \notag \\
P_{D_5} &= \bm{6a_1a_3-12a_4}+u(4a_1a_2-14a_3)-4u^2a_2  \notag \\
& \cdots \;\; \cdots \notag
\end{align}

Return back to the setting with $\varphi=p\circ \iota: M \hookrightarrow W \to N$. 
Note that $V, I$ are naturally extended to $M$: we reset $V$ to be the relative tangent bundle $\ker dp$ and $I$ to be the normal bundle of $M$ in $W$. 
Then $L$ is the zero locus of a section of $\Hom(V, I)$ on $M$ given by $V \subset TW \to I$, and it is now chosen to be generic. 
Let $i: L \hookrightarrow M$ be the inclusion, then $[L]=i_*(1)=c_{top}(V^*\otimes I) \in H^*(M)$, and by the projection formula,  
$[{\eta(M)}]=[\overline{\eta(\varphi)}] =i_*(\,[\overline{\eta(f)}])=i_*(i^*P_\eta(u, a_i))=i_*(1)P_\eta(u, a_i)$. 

\begin{cor} \cite{Kaz03a}, \cite[\S 10]{MPW11}
Given a family of divisors, $\varphi=p\circ \iota: M \hookrightarrow W \to N$, as described above, 
it holds that 
$$\Dual \,[\overline{\eta(\varphi)}] = (u^r-u^{r-1}c_1+\cdots +(-1)^rc_r) P_\eta(u, a_i) \; \in H^*(M)$$
after evaluated by $u= c_1(I)$, $a_i=c_i(V^*\otimes I-V)$ and $c_i=c_i(V)$ with $I=\iota^*TW-TM$ and $V=TW-\varphi^*TN$ and  $r:=\dim V=m-n+1$. 
\end{cor}

For a general Legendre map factored by $L \to PT^*N \to N$, we may take a local chart of $PT^*N$ as $J^1(\C^{n-1}, \C)=T^*\C^{n-1}\oplus \C$ with coordinates $(x, p, z)$ so that  the Legendre immersion is locally described using some generating function $z=g(p)$; the $z$-line $\C$ corresponds to $I$ and the $p$-space $\C^{n-1}_{p}$ corresponds to $V$. We can apply Theorem \ref{thm_tp_L} to this setting \cite[p.731]{Kaz03a}: 

\begin{cor} \cite{Kaz03a}
Given a generic Legendre map $\varphi_L: L \to PT^*N \to N$, 
it holds that 
$$\Dual \,[\overline{\eta(\varphi_L)}] = P_\eta(u, a_i) \; \in H^*(L)$$
after evaluated by $u= c_1(I)$, $a_i=c_i(\varphi_L^*TN- TL - I)$ where $I$ is the dual tautological line bundle over $PT^*N$, equivalently, the normal bundle of the contact hyperplane field.  
\end{cor}

\begin{rem}\upshape  \label{lagrange}
{\bf (Lagrangian case)} 
In the above argument, if one takes $I$ to be the trivial line bundle,  we cover the case of Lagrange maps, which is tied with the $\mathcal{R}$-classification of function-germs. The universal polynomials of type $\eta$ for Lagrange singularities, called {\em Lagrange Thom polynomials},  are obtained by $P_\eta$ with $u=0$. For example, in the above sample, the bold font depicts those parts. 
Lagrange Thom polynomials lie on the polynomial ring generated by $a_1, a_2, \cdots$ with the relation  
$$(1+a_1+a_2+\cdots)(1-a_1+a_2-\cdots)=1$$
which is called the ring of {\em universal Lagrange characteristic classes}. 
It is torsion free and has the additive basis of monomials $a_1^{i_1}a_2^{i_2} \cdots a_s^{i_s}$ with $i_k=0$ or $1$. 
This ring is actually isomorphic to (the limit of) the cohomology ring of the complex Lagrange Grassmannian 
$$\Lambda_n^\C=Sp_n/U_n$$ 
i.e., the space of Lagrange subspaces of the symplectic vector space $\C^{2n}$,  when one replaces $a_i$ by the $i$th Chern class of the tautological bundle of rank $n$ over $\Lambda_n^\C$. 
In Legendrian case, the ring $\mathcal{L}$ is isomorphic to the cohomology ring of the twisted complex Lagrange Grassmannian $\tilde{\Lambda}_n^\C$, that is the space of Lagrange subspaces in the twisted symplectic fibers of the bundle $\C^n \oplus (\C^n \otimes \xi)$, where $\xi$ is the canonical line bundle over $BU_1=\Proj^\infty$. 

Positivity of Lagrange and Legendre Thom polynomials with respect to some particular additive basis (`twisted' Schur $Q$-functions basis and also other basis) has been shown in Mikosz-Pragacz-Weber \cite{MPW09, MPW11}; especially, it contains the list of explicit expansions for $ADE$-types and $P_8, P_9, X_9$ up to codimension $7$. 
There is still room for further studies in this direction. 

\end{rem}

\subsection{Thom polynomials for $\A$-classification}\label{A_classification}
Thom polynomial $tp^\A$ for $\A$-classification of map-germs can also be considered. 
Since for $\A$-equivalence, there is no stabilization of dimension (with respect to suspension), i.e.,  Lemma \ref{stabilization} does not hold for $\A$-orbits $\eta$, thus in general $tp^\A(\eta)$ is just a polynomial in Chern classes of the source bundle and that of the target one, not written in terms of quotient Chern classes. 

A relevant geometric setting for $\A$-classification is described as follows. 
Consider the commutative diagram 

$$\xymatrix{ 
X \ar[rr]^{f} \ar[dr]_{p_0} & & 
Y\ar[dl]^{p_1} 
\\
& B &
}
$$
where $X, Y, B$ are complex manifolds, 
$p_0: X \to B$ and $p_1: Y \to B$ are proper submersions 
of constant relative dimension $m$ and $n$, respectively.  
For each $x \in X$, 
the germ at $x$ of $f$ restricted to the fiber is defined:
$$f|_{p_0^{-1}(p_0(x))}: (\C^m, 0) \to (\C^n,0)$$ 
(in local coordinates of fibers centered at $x$ and $f(x)$, respectively). 
Given an $\A$-finite singularity type $\eta$ of map-germs $(\C^m,0) \to (\C^n,0)$, 
the {\it singularity locus} $\eta(f) \subset X$ 
and the {\it bifurcation locus} $B_\eta(f)=p_0(\eta(f)) \subset B$ 
are defined.

\begin{thm} \cite{HaefligerKosinski, Sasajima18} 
Let $\eta$ be an $\A$-finite singularity type and  $f: X \to Y$ as above. 
If $f$ is appropriately generic, $\Dual [\overline{\eta}(f)] \in H^*(X)$ is expressed 
by  a universal polynomial $tp^{\A}(\eta)$ 
in the Chern class $a_i=c_i(T_{X/B})$ and $b_j=c_j(T_{Y/B})$
of relative tangent bundles. 
$\Dual [\overline{B_\eta(f)]} \in H^*(B)$ is also 
expressed by the pushforward $p_{0*}tp^{\A}(\eta)$. 
\end{thm}
$$\xymatrix{
\overline{\eta}(f) \ar[d]_{p_0} \;\ar@{^{(}->}[r] \; & X \ar[rr]^{jf\qquad \quad} 
\ar[d]_{p_0}  &&  J(T_{X/B}, f^*T_{Y/B})
\\
\overline{B_\eta}(f) \; \ar@{^{(}->}[r] \;& B &&
}
$$
As for Schur basis expansion of $tp^\A$ and the positivity, see \cite{PW08}. 

Let us look at some computational results. 
In case of $m=n=2$, the Thom polynomial of an $\A$-singularity type $\eta$ is defined to be 
$$tp^\A(\eta) \in \Z[a_1, a_2, b_1, b_2]$$
where $a_i, b_i$ are Chern classes of relative tangent bundles of source and target, respectively. 
The $\A$-orbits with $\A$-codimension $\le 4$ are all simple (no nearby moduli). 
For each of fold, cusp, swallowtail, butterfly and $I_{2,2}^{1,1}$, 
the $\A$-orbit is an open dense subset of its $\K$-orbit in $J(2,2)$,  that implies that 
$tp^\A$ coincides with $tp$ for its $\K$-orbit. 
For other $\A$-type $\eta$, the $\A$-orbit has positive codimension in its $\K$-orbit, thus   
$tp^\A(\eta)$ can not be expressed in terms of quotient Chern classes $c_i=c_i(f)$. For instance, $tp^{\A}$ for $\A$-types, lips(beaks) $(x^3+ xy^2,y)$, gulls $(x^4+xy^2+x^5,y)$, butterfly $(x^5+ xy+x^7,y)$, 
are computed by Sasajima-Ohmoto  \cite{Sasajima18} with application to classical projective algebraic geometry:  

\begin{eqnarray*}
tp^\A(\mbox{lips}) &=&
{-2a_1^3+5a_1^2b_1-4a_1b_1^2-a_1a_2+a_2b_1+b_1^3}\\
tp^\A(\mbox{gulls}) &=& 
{6a_1^4 -a_1^2a_2-4a_2^2-17a_1^3b_1+4a_1a_2b_1+17a_1^2b_1^2-3a_2b_1^2} \\
&&
{-7a_1b_1^3+b_1^4+2a_1^2b_2+6a_2b_2-4a_1b_1b_2+2b_1^2b_2-2b_2^2 } \\
 tp^\A(\mbox{goose}) &=&   
{2a_1^4 +5a_1^2a_2+4a_2^2-7a_1^3b_1-10a_1a_2b_1+9a_1^2b_1^2+5a_2b_1^2} \\
&&
{-5a_1b_1^3+b_1^4-2a_1^2b_2-6a_2b_2+4a_1b_1b_2-2b_1^2b_2+2b_2^2 } \\
\end{eqnarray*}

Not only for plane-to-plane germs, the $\A$-Thom polynomials $tp^\A$ are also computed in \cite{Sasajima18} for other $\A$-simple-classifications of $(\C^2,0) \to  (\C^3,0)$, $(\C^3,0) \to  (\C^3,0)$, and $(\C^3,0) \to (\C^4,0)$. 
There have been $\A$-simple classifications of parametric curves in plane, $3$-space, etc., and those are collectively united by Arnold into a stable form of $\A$-classification of germs $(\C,0) \to (\C^n,0)$ with involving $n$ as a parameter, and the corresponding Thom polynomial theory would be a new one, see Pissolato \cite{Pissolato}. 

\begin{rem}\upshape
{\bf (Hurwitz space)} 
The case of maps between families of curves has intensively been explored by Kazarian-Lando \cite{KL04, KL07} and Kazarian-Lando-Zvonkine \cite{KLZ18, KLZ22} for the study of Hurwitz spaces and genus $0$ Gromov-Witten theory. Those are actually more advanced ones than ordinary Thom polynomial theory, because the family are allowed to have nodal degeneration and orbifold structures are involved. It also suggests Thom polynomials corresponding to certain classification of function/map-germs over isolated complete intersection singularities.
\end{rem}

\section{Thom polynomials for multi-singularities}\label{multising_tp}

At the beginning of this century, another breakthrough occurred. Based on his topological argument using complex cobordism theory $MU^*$, M. Kazarian \cite{Kaz03b, Kaz06} proposed a general framework for  universal polynomial expression of the multi-singularity locus of prescribed type for proper holomorphic maps and showed its high potential in application. 
We briefly describe this conjectural theory and report a progress towards an entire proof.

\subsection{Multi-singularity loci and Hilbert schemes}\label{multising}
Below we use the terminology `multi-singularity' in the following sense: 
\begin{definition}{\rm 
A {\it  multi-singularity type of map-germs with codimension $\kappa$} is an ordered set 
$\underline{\eta}:=(\eta_1, \cdots, \eta_r)$ 
of $\K$-types $\eta_i$ of mono-germs $(\C^m, 0) \to (\C^n,0)$ for some $m, n$ with $\kappa=n-m$. 
In case of $\kappa\le 0$, we assume that 
the collection  $\underline{\eta}$ contains no submersion-germs. 
}
\end{definition}
Note that the entries $\eta_i$ can be duplicate. 
We denote by $\Aut(\eta_1, \cdots, \eta_r)$ the subgroup of permutations in $\mathfrak{S}_r$ preserving the types $\eta_i$, e.g., $\Aut(A_1, A_1, A_2, A_1 ,A_2)\simeq \mathfrak{S}_3 \times \mathfrak{S}_2$.  

A multi-jet bundle ${}_rJ^k(M, N)$ is defined by the restriction of the cartesian product of $r$ copies of $J^k(M, N) \to M \times N$ over  $F(M, r) \times N^r$,  where $F(M, r)=M^r \smallsetminus \Delta^r_M$ is the configuration space of ordered distinct $r$ points of $M$ (i.e., the complement of the fat diagonal $\Delta^r_M$) \cite{AGV, MondNuno}.  Put $\pi^r_M: {}_rJ^k(M, N) \to F(M, r)$ and $\pi^r_N: {}_rJ^k(M, N) \to N^r$ to be the canonical projections. For a map $f: M \to N$,  the multi-jet extension map ${}_rj^kf$ is defined by 
$${}_rj^kf: F(M, r) \to {}_rJ^k(M, N), \quad (p_1, \cdots, p_r) \mapsto (j^kf(p_1), \cdots, j^kf(p_r)).$$  

Let $\underline{\eta}:=(\eta_1, \cdots, \eta_r)$ be a multi-singularity type, and $\dim M=m$ and $\dim N=n \, (=m+\kappa)$. 
We always assume that $k$ is sufficiently high so that it is greater than the determinacy degree of every entry $\eta_i$. 
Let $W_{\underline{\eta}} \subset {}_rJ^k(M, N)$ be the submanifold consisting of $z=(z_1, \cdots, z_r)$ such that $z_i$ is the $k$-jet of a mono-germ of type $\eta_i$ ($1\le i \le r$) and $\pi_N^r(z)$ lies on the diagonal $\Delta^r_N$ of $N^r$. 
Let $S=\{p_1, \cdots, p_r\} \subset M$, a set of $r$ distinct points, and suppose that $f(p_1)=\cdots = f(p_r)=:q$ and the germ $f: (M, p_i) \to (N, q)$ is of type $\eta_i$ for each $1 \le i \le r$. 
The multi-germ $f: (M, S) \to (N, q)$ is {\em stable} if the multi-jet extension ${}_rj^kf$ (of a representative of the germ $f$) is transversal to $W_{\underline{\eta}}$ at a point $(p_1, \cdots, p_r) \in F(M, r)$.

\begin{definition}\upshape
We say that $f: M \to N$ is a {\em locally stable map} (or simply, a stable map) if for any finite subset $S$ of $M$ such that $f$ maps $S$ to a point $q \in N$, the germ $f: (M, S) \to (N, q)$ is stable. 
\end{definition}

Let $f: M \to N$ be a proper and locally stable map. 
Take the preimage of $W_{\underline{\eta}}$ via ${}_rj^kf$ and project it to the first factor $M$ of $F(M, r) \subset M^r$, then the resulting set is written as  
$$\underline{\eta}(f)
:=\left\{ \; p_1 \in \eta_1(f) \; \Bigl| \;
\begin{array}{l} 
\exists\,  p_2, \cdots, p_r \in f^{-1}f(p_1)\smallsetminus \{p_1\} \; s.t.  \; p_i\not=p_j\; (i\not=j) \\
 \mbox{and the germ of $f$ at $p_i$ is of type $\eta_i$ } 
\end{array}  \right\},$$
which becomes a locally closed submanifold of $M$ having the expected dimension. 
This locus is mapped via $f$ to a locally closed submanifold of $N$ of the same dimension: 
$$f(\underline{\eta}(f))
:=\left\{ \; q \in N \; \Bigl|  \;
\begin{array}{l} 
\exists\, p_1, \cdots, p_r \in f^{-1}(q)  \; s.t.  \; p_i\not=p_j\; (i\not=j) \\
 \mbox{and the germ of $f$ at $p_i$ is of type $\eta_i$} 
\end{array}  \right\}.$$
Note that $f: {\underline{\eta}(f)} \to {f(\underline{\eta}(f)})$ 
is a finite-to-one covering whose degree, denoted by 
$\deg_1 \underline{\eta}$,  is equal to the number of $\eta_1$ appearing in the tuple $\underline{\eta}$. 

For our convenience, we introduce the {\em source and target multi-singularity loci classes of $f$} to be 
\begin{equation}\label{multi}
m_{\underline{\eta}}(f):= \#\Aut(\eta_2, \cdots, \eta_r) [\overline{\underline{\eta}(f)}], \;\; 
n_{\underline{\eta}}(f):= \#\Aut(\eta_1, \cdots, \eta_r) [\overline{f(\underline{\eta}(f))}] 
\end{equation}
(more generally, these classes can be defined for {\em arbitrary} proper maps  \cite{Ohmoto24}). 
In particular, it holds that 
$$f_*m_{\underline{\eta}}(f) = n_{\underline{\eta}}(f).$$

Based on his topological observations supported by many computations, Kazarian claimed in \cite{Kaz03b, Kaz06} the following statement:

\

\t
{\rm \bf Thom-Kazarian principle for multi-singularities:} 
{\em 
For any  multi-singularity type $\underline{\eta}=(\eta_1, \cdots, \eta_r)$ of map-germs of relative codimension $\kappa$, there is a unique abstract Chern polynomial 
$$R_{\underline{\eta}}=\sum a_I c^I \in \Q[c_1, c_2, \cdots],$$ 
called the residual polynomial of type $\underline{\eta}$, 
such that $R_{\underline{\eta}}$ does not depend on the order of entries in $\underline{\eta}$ 
and recursively admits the following property: 
For a proper morphism $f: M \to N$ with $\dim N-\dim M=\kappa$, 
the  $\underline{\eta}$-singularity loci classes $m_{\underline{\eta}}(f)$ and $n_{\underline{\eta}}(f)$ 
are expressed by 
$$m_{\underline{\eta}}(f) = \sum R_{J_1}\cdot f^*f_*(R_{J_2}) \cdots f^*f_*(R_{J_s}) \eqno{(stp)}$$
$$n_{\underline{\eta}}(f) =\sum f_*(R_{J_1})\cdots f_*(R_{J_s})\eqno{(ttp)}$$
in rational cohomology rings of $M$ and $N$, 
where the summand runs over all possible partitions of the set $\{1, \cdots, r\}$ 
into a disjoint union of non-empty unordered subsets, 
$\{1, \cdots, r\}= J_1\sqcup \cdots \sqcup J_s \; (s\ge 1)$, 
and the subset containing the entry $1 \in \{1, \cdots, r\}$ is denoted by $J_1$. 
Here, for each $J=\{j_1, \cdots, j_k\}$,  
$R_J$ stands for the residual polynomial $R_{(\eta_{j_1}, \cdots, \eta_{j_k})}$ 
evaluated by $c_i=c_i(f)$, the quotient Chern classes associated to $f$. 
}

\

We call the universal polynomials in the right hand side of the formulas $(stp)$ and $(ttp)$ 
the {\it source and target Thom polynomial of a stable multi-singularity type $\underline{\eta}$}, respectively. 
Note that the source Thom polynomial $m_{\underline{\eta}}(f)$ is written in quotient Chern classes $c_i(f)$'s and Landweber-Novikov classes $f^*s_I(f)$'s, and the target one  is in $s_I(f)$'s, 
and those polynomials depend only the multi-singularity type $\underline{\eta}$. 
For convenience, we introduce a reduced version of multi-singularity Thom polynomials 
(the last equality is for locally stable maps $f$): 
$$ tp(\underline{\eta})(c(f)):=\frac{1}{\#\Aut(\eta_2, \cdots, \eta_r)}m_{\underline{\eta}}(f)=[\overline{\underline{\eta}(f)}],$$ 
$$ tp_{\rm target}(\underline{\eta})(c(f)):=\frac{1}{\#\Aut(\eta_1, \cdots, \eta_r)}n_{\underline{\eta}}(f)=[\overline{f(\underline{\eta}(f))}].$$ 
Note that 
$tp_{\rm target}(\underline{\eta})(c(f)):=\frac{1}{\deg_1\underline{\eta}}f_*tp(\underline{\eta})(c(f))$ where $\deg_1\underline{\eta}$ is the number of entries in $\underline{\eta}$ which is the same type as $\eta_1$. 

\begin{exam}\label{exam33} 
\upshape
Here are some samples of Thom polynomials for multi-singularities of stable maps in low dimensions \cite{Kaz03b, Kaz06, Ohmoto16, RimanyiPortal}. Using the Thom-Kazarian principle, Rim\'anyi's interpolation method works effectively  
for computing residual polynomials $R_{\underline{\eta}}$. 
First, in case of $\kappa=0$, we see 
\begin{align*}
tp(A_1)  &= c_1 \\
tp(A_2)  &= c_1^2 + c_2 \\
tp(A_1^2)  &= c_1s_1-4c_1^2-2c_2\\
tp(A_3)  &= c_1^3+3c_1c_2+2c_3\\
tp(A_1^3)  &= \textstyle \frac{1}{2}(c_1s_1^2-4c_2s_1-4c_1s_2-2c_1s_{01}-8c_1^2s_1+40c_1^3+56c_1c_2+24c_3)\\
tp(A_1A_2)  &=  c_1s_2+c_1s_{01}-6c_1^3-12c_1c_2-6c_3\\
tp(A_2A_1)  &= c_1^2s_1+c_2s_1 -6c_1^3-12c_1c_2-6c_3.
\end{align*}
Also in case of $\kappa=1$, 
\begin{align*}
tp(A_0^2)  &= s_0-c_1 \\
tp(A_1)  &= c_2 \\
tp(A_0^3)  &= \textstyle \frac{1}{2}(s_{0}^2-s_1-2s_0c_1+2c_1^2+2c_2)\\
tp(A_0A_1)  &= s_{01}-2c_1c_2-2c_3\\
tp(A_1A_0)  &=  s_0c_2-2c_1c_2-2c_3 \\
tp(A_0^4)  &= \textstyle  
\frac{1}{6}(s_0^3-3s_0s_1+2s_2+2s_{01}-3s_0^2c_1+3s_1c_1+6s_0c_1^2+6s_0c_2-6c_1^3\\
&\quad -18c_1c_2-12c_3).
\end{align*}
\end{exam}

Regarding the validity of the above principles, there are some evidences, attempts and solid results, some of which are ongoing projects, as mentioned below. 
Commonly among those different approaches, it is important to find a nice compactification of the configuration space $F(M, r)$ and to rigorously define multi-singularity loci classes in it --  the Hilbert schemes $M^{[l]}$ of $0$-dimensional subschemes with colength $l$ of $M$ takes an essential role here.

\

\t
$\bullet\;\;$ 
{\em Double and triple-point formulas} (Ronga \cite{Ronga73, Ronga84}, Fulton-Laksov \cite{FultonLaksov}): \\
For {\em arbitrary} proper maps $f: M \to N$ with $\kappa >0$, the {\em double and triple-point loci classes $m_{A_0^2}(f)$ and $m_{A_0^3}(f)$ }are defined, and it is shown that 
$$m_{A_0^2}(f)=f^*f_*(1)-c_{\kappa}(f),$$
$$m_{A_0^3}(f)=f^*f_*(m_{A_0^2}(f))-2c_\kappa(f) f^*f_*(1)+2\left(c_\kappa(f)^2+\sum_{i=0}^{\kappa-1}2^ic_{\kappa-i-1}(f)c_{\kappa+i+1}(f)\right).$$
These are re-written in the form of $(stp)$ using Landweber-Novikov classes. 
The proofs of the double and triple-point formulas rely on the desingularization method. Namely, desingularizations of the loci are explicitly constructed by blowing-ups along diagonals, and actually the resulting spaces coincide with relative Hilbert schemes associated to $f$, which are defined in the ordered Hilbert schemes $M^{[[l]]}$ for $l=2, 3$ (cf.  \cite{Gaffney, Ohmoto24}). Following this line, the {\em quadruple-point formula} was pursued by Dias-Le Barz \cite{DL}, but it meets a serious difficulty caused from 
the fact that the Hilbert scheme $M^{[[4]]}$ is singular. 
On the other hand, once we assume the Thom-Kazarian principle, the quadruple-point formula is explicitly computed by Marangell-Rim\'anyi \cite{MR10} via the interpolation.  \\

\t
$\bullet\;\;$ 
{\em Curvilinear multiple-point theory} (Kleiman \cite{Kleiman81, Kleiman90}, Colley \cite{Colley}, Ran \cite{Ran85}): \\
The principle is strongly supported by the multiple-point theory for an important particular class of  {\it curvilinear maps}, that is, maps $f: M \to N$  with $\kappa \ge 0$ and corank at most one; its singularity type is denoted by $A_\mu$, i.e., the local ring is isomorphic to $k[x]/\langle x^{\mu+1}\rangle$, and an $r$-tuple of them is denoted by $A_\bullet=(A_{\mu_1}, \cdots, A_{\mu_r})$.  That has been established by Kleiman and others around the 80s, where 
main tools are curvilinear Hilbert schemes and the iteration method \cite{Kleiman81, Kleiman90}. 
In this case,  moduli spaces under consideration become non-singular, and thus concrete computations could be done. 
For instance, there is a recursive formula for $m_{A_0^r}$ (multiple-point formulas) by Kleiman's iteration method  and 
Colley's program produces $m_{A_\bullet}$ of low codimensions. 
Note that those formulas are only applicable to curvilinear maps,  i.e., they give wrong answers for non-curvilinear maps. 
\\

\t
$\bullet\;\;$ 
{\em Tautological integrals over curvilinear Hilbert schemes} (B\'erczi-Szenes \cite{BercziSzenes21}): \\
The target multiple-point formula for generic maps (including non-curvilinear maps) has been established. That corresponds to the formula $(ttp)$ for the type $\underline{\eta}=A_0^r$ where $f$ is a proper locally stable map. Their method relies on torus equivariant localization technique -- they use iterated residue integrals with partial desingularization based on the curvilinear Hilbert schemes and non-reductive group quotient $J^k(1, n)\git \Aut(\C^n, 0)$, that has been developed as a byproduct from their main project for computing Thom polynomials for $A_\mu$-singualrities (Theorem \ref{thm_BS}).  \\

\t
$\bullet\;\;$ 
{\em Functorial approach using algebraic cobordism} (Ohmoto \cite{Ohmoto24}): \\
In contrast to the approaches above, this is a non-constructive approach which proves the existence of universal expression $(ttp)$, and partly $(stp)$, in the Thom-Kazarian principle. In particular, the results ensure the correctness of many formulas counting multi-singular points which are obtained via the interpolation under assuming the principle, including the quadruple-point formula of Marangell-Rim\'anyi \cite{MR10} for generic maps. 
Given a multi-singularity type $\underline{\eta}=(\eta_1, \cdots, \eta_r)$,  the {\em source and target multi-singularity loci classes} $m_{\underline{\eta}}(f)$ and $n_{\underline{\eta}}(f)$ are defined for {\em arbitrary} proper maps $f: M \to N$ by using intersection theory on relative Hilbert schemes so that they coincide with the previously defined ones in (\ref{multi}) if $f$ is locally stable. Then, the existence theorem of target multi-singularity Thom polynomials $(ttp)$ for the loci class $n_{\underline{\eta}}(f)$ has been proven, by developing an idea of Kazarian \cite{Kaz03b} in the context of algebraic geometry -- a main feature of the proof is a striking use of algebro-geometric {\em cohomology operations}, established by A. Vishik \cite{Vishik} (the universality of Landweber-Novikov operations), on algebraic cobordism theory $\Omega^*$ of Levine-Morel \cite{LevineMorel, LevinePandhari}. 
Namely, it turns out that 
$$n_{\underline{\eta}}: \Omega^\kappa(N) \to \CH^{\ell}(N), \quad [f: M \to N] \mapsto n_{\underline{\eta}}(f)$$
is a well-defined cohomology operation, and thus, it is uniquely written as a polynomial of $s_I(f)$'s after tensoring by $\Q$. 
 On the other hand, due to some technical difficulty,  the existence of the source multi-singularity Thom polynomials  $(stp)$ has been shown, up to now, only in a certainly nice case, i.e., for proper locally stable maps $f: M \to N$ with $\kappa \ge 0$ plus some additional conditions -- the proof relies on  the existence theorem $(ttp)$ and the interpolation method of Rim\'anyi. 
For finding an entire proof of the Thom-Kazarian principle for $m_{\underline{\eta}}(f)$, most critical is that the classifying space of multi-singularities has not yet been found properly. 
That suggests, following Kazarian's view  \cite{Kaz03b}, to build a bridge between Thom polynomial theory and motivic homotopy theory.

\subsection{Higher Thom polynomials for multi-singularity types}\label{higher_multising}
It is natural to expect the Thom-Kazarian principle also holds for the SSM class of the multi-singularity locus of prescribed type, that is, we conjecture that there is a universal polynomial (series) $tp^{\SM}(\overline{\underline{\eta}})$ in variables $c_i$'s and $s_I$'s satisfying 
$$s^{\SM}(\jeden_{\overline{\underline{\eta}(f)}}, M) = tp^{\SM}(\overline{\underline{\eta}})(c(f))$$
in $H^*(M; \Q)$ evaluated by $c_i=c_i(f)$ and $s_I=f^*s_I(f)$ for any proper locally stable maps $f$. 
The leading term of the SSM class is the fundamental class $[\overline{\underline{\eta}(f)}]$, thus it is expressed 
by the reduced multi-singularity Thom polynomial, namely, 
$$tp^{\SM}(\overline{\underline{\eta}})(c(f))=tp(\underline{\eta})(c(f)) + h.o.t.$$
By definition, $c(TM)\cdot  tp^{\SM}(\overline{\underline{\eta}})(c(f))$ coincides with the CSM class $C_*(\jeden_{\overline{\underline{\eta}(f)}})$. 
Multiply the both sides of the equality by $\bar{c}(f)=c(f)^{-1}=c(TM-f^*TN)$ and then apply the pushforward $f_*$, 
we see
$$s^{\SM}(f_*(\jeden_{\overline{\underline{\eta}(f)}}), M) = f_*( \bar{c}(f) \cdot tp^{\SM}(\overline{\underline{\eta}})(c(f)))$$
by using the naturality $f_*C_*=C_*f_*$ of the CSM class transformation, the definition of SSM class and the projection formula. 
Notice that this total cohomology class in the target space $N$ is written in terms of the Landweber-Novikov classes $s_I$. 

Recall the adjacency relation of multi-singularities both in source and target: 
the arrow $\underline{\eta} \to \underline{\xi}$ means that the stratum of type $\underline{\xi}$ is contained in the closure of the stratum of $\underline{\eta}$. That yields the set of all multi-singularity types to be a poset (partially ordered set), and we call it the {\em multi-singularity poset}. Note that this makes sense for each of the source and the target spaces according to the definition of loci $\underline{\eta}(f)$ and $f(\underline{\eta}(f))$ for locally stable maps $f$; in the source space, the first entry $\eta_1$ is distinguished from the others, and in the target space, the order of all entries is ignorable. 

For instance, the Hasse diagram of the source multi-singularity poset in case of $\kappa=0$ starts as follows: 

{\small 
$$\xymatrix{
A_1 \ar[r]  \ar[rd] & A_1^2 \ar[r] \ar[rd] \ar[rdd] &A_1^3  \ar[r]  \ar[rd] \ar[rdd] &  \cdots\\
& A_2 \ar[r]  \ar[rd] & A_2A_1\, \&\, A_1A_2\ar[r] \ar[rd] & \cdots \\
&& A_3 \ar[r]   & \cdots \\
}
$$
}

Assuming the Thom-Kazarian principle, the interpolation method of Rim\'anyi yields several computations
of the SSM classes for multi-singularity loci provided $\kappa$ is fixed. 
Furthermore, in a sophisticated manner inspired by geometric representation theory and the previous work \cite{FR18, Rimanyi25},  Thom series expression of the SSM classes with parameter $\kappa$ is also studied in a recent paper of Koncki-Rim\'anyi \cite{KR25}. 

For instance, in case of $\kappa=0$, some sample computations are given as follows \cite{Ohmoto16, Nekarda}: 
\begin{align*}
tp^{\SM}(\overline{A_1})
&= \textstyle c_1-c_1^2+c_1^3+\cdots,\\
tp^{\SM}(\overline{A_2})
&= \textstyle c_1^2+c_2-2c_1^3-3c_1c_2-c_3+\cdots,\\
tp^{\SM}(\overline{A_1^2})
&= \textstyle c_1 s_1 -4 c_1^2 - 2 c_2-4 c_1^3 - 10 c_1 c_2 - 4 c_3 + c_1 s_{01} \\
&\quad \textstyle
+3 c_1^2 s_1 + 2 c_2 s_1 - \frac{1}{2}c_1 s_1^2+\cdots,\\
tp^\SM(\overline{A_3})
&=c_1^3+3c_1c_2+2c_3+\cdots, \\
tp^\SM(\overline{A_2A_1})
&=s_1c_1^2-6c_1^3+s_1c_2-12c_1c_2-6c_3+\cdots, \\
tp^\SM(\overline{A_1A_2})
&=s_2c_1+s_{01}c_1-6c_1^3-12c_1c_2-6c_3+\cdots,\\
tp^\SM(\overline{A_1^3}) 
&=\textstyle \frac{1}{2}(s_1^2c_1-4s_2c_1-2s_{01}c_1-8s_1c_1^2 -4s_1c_2\\
&\quad \textstyle + 40c_1^3 + 56c_1c_2 + 24c_3)+\cdots
\end{align*}

For instance, assume that $\kappa \le 0$ and 
let $D(f) := \overline{f(A_1(f))}$ be the discriminant locus.  
We call $C_*(\jeden_{D(f)})\in H^*(N)$ {\it the discriminant Chern class of $f$}. 
By using the M\"obius inversion formula associated to the multi-singularity poset, 
we find a unique (universal) constructible function $\disconst  \in \F(M)\otimes \Q$ 
so that $f_*\disconst = \jeden_{D(f)}$. 
For instance, in case of $\kappa=0$, we compute it as 
$$\textstyle \disconst
:=\jeden_{\overline{A_1}}-\frac{1}{2}\jeden_{\overline{A_1^2}}
-\frac{1}{6}\jeden_{\overline{A_1^3}}+\frac{1}{2}\jeden_{\overline{A_3}} 
+ \cdots \;\; \in \F(M)\otimes \Q.$$

The Thom-Kazarian principle indicates that 
there should be a polynomial $tp^{\SM}(\disconst)$ 
in variables $c_i=c_i(f)$ and $s_I=s_I(f)$ such that 
$$
C_*(\disconst)=c(TM)\cdot tp^{\SM}(\disconst) \;\; \in H^*(M; \Q)
$$
and also a polynomial $tp^{\SM}(\jeden_{D(f)})$ in $s_I=s_I(f)$ such that 
$$ 
C_*( \jeden_{D(f)})=c(TN)\cdot tp^{\SM}( \jeden_{D(f)})  \;\; \in H^*(N; \Q). 
$$
In particular, the Euler characteristics of the discriminant is universally expressed by 
$$
\chi(D(f))
=\int_M c(TM) \cdot  tp^{\SM}(\disconst)
=\int_N c(TN)\cdot tp^{\SM}(\jeden_{D(f)}).$$
The naturality of $C_*$ appeals to the inclusion-exclusion principle (indeed, a primitive form of such an argument was often seen in classical computations of elementary characters of projective varieties). 
For example, in case of $\kappa=0$,  it follows from Table \ref{table_tpSM} and the above expression of $\disconst$ 
that the low degree terms of the SSM class are given by 
\begin{eqnarray*}
tp^{\SM}(\disconst)
&=&{\textstyle c_1+\frac{1}{6}(6c_1^2+6c_2-3c_1s_1)}\\
&&{\textstyle +\frac{1}{6}\left(
\begin{array}{l}
{\textstyle c_1^3+11c_1c_2+6c_3-2c_1s_{01}-5c_1^2 s_1}\\
{\textstyle  -4 c_2s_1+c_1s_1^2+2c_1s_2 }
\end{array}
\right) + h.o.t.}
\end{eqnarray*}

A reduced divisor $D$ in a complex manifold $N$ is called to be 
{\it free} (in the sense of Kyoji Saito) 
if the sheaf of germs of logarithmic vector fields ${\rm Der}_N(-\log D)$ 
is locally free. 
As for the CSM class of a free divisor $D$,  
the following equality was conjectured  by P. Aluffi, 
and was proved by X. Liao \cite{Liao}: 

\begin{thm}[{\bf CSM class of free divisors} \cite{Liao}] 
If $D$ is locally quasi-homogeneous, it holds that 
\begin{center}
$c^{\SM}(N-D)=c({\rm Der}_N(-\log D))  \quad  \in H^*(N).$
\end{center}
\end{thm}

It is well-known that for a stable germ with $\kappa \le 0$, i.e., universal unfolding of ICIS germs, the discriminant  is a free divisor (cf. \cite{MondNuno}). 
On the other hand, we have seen that 
the CSM class of $D=D(f)$ in the ambient space $N$ is expressed  
using the target universal SSM class for $f_*\disconst = \jeden_{D(f)}$: 
\begin{eqnarray*}
c^{\SM}(N-D)=c^{\SM}(N)-c^{\SM}(D(f))
=c(TN)(1-tp^{\SM}(\jeden_{D(f)})). 
\end{eqnarray*}
Thus we get the following formula: 
$$tp^{\SM}(\jeden_{D(f)})=1-\frac{c({\rm Der}_N(-\log D(f))}{c(TN)}.$$
This is an unexpected and highly non-trivial identification of the target SSM Thom polynomial $tp^{\SM}(\jeden_{D(f)})$ in variables $s_I=s_I(f)$, which measures the difference of Chern classes for locally free sheaves of vector fields and of logarithmic vector fields. That might be comparable to the conjectural formula (\ref{Milnor1}) on the Milnor number constructible function.

\section{Some topics in applications}\label{applications}
\subsection{Application to enumerative geometry -- instanton counting}\label{enumerative_geometry}
The literature on enumerative geometry contain full of counting formulas, including 19th century classics as well as modern ones inspired by mathematical physics and applied algebraic geometry. 
Thom polynomial theory serves as a kind of {\em taxonomy} or {\em systematics} to re-organize these counting formulas and discover new species. 

Here we present a typical example in this direction -- the enumeration of singular curves on surfaces -- according to Kazarian \cite[\S 10]{Kaz03b}.  This classical topic has been resurrected since the 90s due to new interest from string theory, known as {\em instanton counting}. 
Consider the linear system $|\mathcal{L}|=\Proj(H^0(S;\mathcal{L}))$ associated to a line bundle $\mathcal{L}$ on a projective non-singular complex surface $S$. 
Let $\underline{\eta}=(\eta_1, \cdots, \eta_r)$ be a combination of types of plane curve singularities $(\C^2, 0) \to (\C, 0)$, e.g., the type of $r$ nodes denoted by $A_1^r=(A_1, \cdots, A_1)$, and put $n$ to be the expected codimension of the target multi-singularity locus of type $\underline{\eta}$ (e.g., it is $r$ for $A_1^r$ with $\kappa<0$).  
Then, a general problem is to determine the number of singular curves in this system which have singular points of the prescribed type $\underline{\eta}$ and pass through $\dim |\mathcal{L}|-n$ common points of $S$ in general position. 

The problem is interpreted as counting $\underline{\eta}$-type multi-singular points of a certain map.  Let $N = \Proj^n$ be a generic linear subspace of $|\mathcal{L}|$ which is transverse to the discriminant. Take the trivial fibration $p: N\times S \to N$ and the incident hypersurface $\iota: M \hookrightarrow N\times S$ with the normal bundle $I$, i.e., $M$ consists of the pairs of a curve in the family and a point lying on this curve. Now $\dim M=n+1$ and $\dim N=n$ (thus $\kappa=-1$). 
Let $f=p \circ \iota: M \to N$ be the composed map, that is the family of curves over the sub-linear system $N$, and put $L=\Sigma^2(f)$,  the critical locus of $f$ (recall the setup in \S \ref{legendre}). 

Now we can apply target Thom polynomials for $n_{\underline{\eta}}$ to this map $f: M \to N$, since the existence theorem of the universal polynomials has been established in \cite{Ohmoto24}. 
In this geometric setting, as remarked in \S \ref{legendre}, the quotient Chern class $c(f)$ is rewritten in terms of Legendre characteristic classes in $\mathcal{L}$, so we state the computation in this terminology following \cite{Kaz03b}. 

Let $h=c_1(\mathcal{L}) \in H^2(S)$, $1+c_1 +c_2 = c(T^*S)$, and let $t$ be the divisor class in $N=\Proj^\ell$. Then
$$u=c_1(I)=[M]=t+h, \;\; a=c(T^*S\otimes I-TS)=\frac{(1+u)^2+(1+u)c_1+c_2}{1-c_1 +c_2}.$$
The pushforward $f_*:H^*(L) \to H^*(N)$ is factorized by $f_* = p_*i_*$, where $i : L \to N\times S$ is the embedding and $p : N\times S \to N$ is the projection to the first factor. 
The class $[L]=i_*(1)$ is easily computed, and the homomorphism $p_*$ commutes with multiplication by $t$ and sends monomials of degree $2$ in the classes $h, c_1, c_2$ to the numbers; 
$d := p_*(h^2) = \mathcal{L} \cdot  \mathcal{L}$,  
$k := p_*(hc_1) = \mathcal{L} \cdot K_S$, 
$s:=p_*(c^2_1)=K_S \cdot K_S$ and 
$x:=p_*(c_2)=\chi(S)$.

Since the target singularity locus of type $\underline{\eta}$ has codiemsion $n$,  
there are finitely many singular curves having isolated singular points of type $\underline{\eta}$ which belong to the family $f$. 
Applying the established Thom-Kazarian principle, the degree $\int n_{\underline{\eta}}(f)$ expresses the number of those singular curves; especially, for any $J \subset \underline{\eta}$, the residual class $S_{J}=f_*R_{J} \in H^*(\Proj^n)$ is a linear combination of the numbers $d, k, s, x$ multiplied by $t^{n_J}$ ($n_J$ is the expected codimension of the locus of type $J$), thus the desired degree is written by a polynomial in  $d, k, s, x$. That is summarized as follows.

\begin{thm}  \label{curve_counting} \cite{Kaz03b} 
Let $\underline{\eta}=(\eta_1, \cdots, \eta_r)$ be a multi-singularity type of plane curve-germs with $\K_e$-codimension $n$. Then it holds that for any projective surface with an ample line bundle $\mathcal{L}$ and any sufficiently generic linear sub-system $\Proj^n \subset |\mathcal{L}|$,  the number $N_{\underline{\eta}}$ of singular curves belonging to the system that  possess the combination of singularities of type $\underline{\eta}$ is given by 
$$N_{\underline{\eta}}=\frac{1}{\#\Aut(\underline{\eta})} \sum \int_{\Proj^n} S_{J_1}\cdots S_{J_s}$$
where the sum runs over all partition of $\{1, \cdots, r\}$ and $S_{J_s}$ denotes the image of the residual polynomial of multi-singularities $\{\eta_j\, | \, j \in J_s\}$. 
\end{thm}

A generating function of the numbers $N_{\underline{\eta}}$ is also formulated \cite{Kaz03b}. 
Theorem \ref{curve_counting} for the simplest type  $\underline{\eta}=A_1^r$ 
solves the problem on the existence of universal polynomials for counting nodal curves, known as part of  the {\em G\"ottsche conjecture} \cite{Goettsche, KP99, KP04} (now several different proofs are available, e.g., \cite{LT14, Rennemo, Berczi17, Berczi23}). 
Like as this, there are a lot of questions for which Thom polynomial theory may effectively give satisfactory answers. Other examples are: 
\begin{itemize}
\item[-]  higher dimensional analogy to G\"ottsche conjecture, i.e., counting divisors with a prescribed combination of singularities in a given linear system;  
\item[-]  other type of contact problems, e.g.,  
 enumerating multi-tangent hyperplanes to a given hypersurface in projective space etc.; 
\item[-]  computing the degree of bifurcation sets, e.g., event surfaces and its boundary strata when viewing a curve and a surface in $3$-space (partly related to Computer Vision); 
\end{itemize}
and so on, see, e.g.,  Vainsencher \cite{Vainsencher, Vainsencher07}, Tzeng \cite{Tzeng12, Tzeng19}, Rennemo \cite{Rennemo},  Kazarian \cite{Kaz03b}, Wall \cite{Wall09}, Nekarda \cite{Nekarda, NekardaOhmoto} and references therein.

\subsection{Vanishing topology of map-germs}\label{vanishing_top}
As a well-known theorem due to J. Milnor, the Milnor fiber of an isolated hypersurface-singularity $f: (\C^{d+1},0) \to (\C,0)$ is homotopic to a bouquet of $d$-dimensional real spheres; the number $\mu(f)$ of the spheres is called the {\em Milnor number} of the germ $f$. 
The number is interpreted as the number of $A_1$ (Morse)-singular points appearing in a stable perturbation of the germ $f$, and also it is computed by the dimension of the Milnor algebra: 
$$\mu(f)= \#  A_1(f_t) = \dim_\C \Ost_{\C^m, 0}/J_f.$$

This picture is generalized to the case of map-germs. 
For the simplicity, we assume that $m \le n$ throughout this subsection.  
Let $f: (\C^m,0) \to (\C^n,0)$ be a $\mathcal{G}$-finitely determined map-germ, and $\underline{\eta}$ a stable (mono/multi-)singularity type  of the expected codimension $n$ in the target. 
Suppose that there is a {\it stable perturbation} 
$$f_t: U \to \C^n\;\; (t\in \Delta \subset \C, \; 0 \in U \subset \C^m)$$
so that $f_0$ is a  representative of the germ $f$ and $f_t$ for $t\not=0$ is a stable map. 
Note that if $\eta$ is a mono-singularity type,  it is enough to assume that $f$ is $\K$-finitely determined, while for a multi-singularity type, we need the $\A$-finiteness of $f$. Then $\underline{\eta}(f_t)$ for $t\not=0$ consists of finitely many isolated points; 
the number is constant for non-zero $t$ and does not depend on the choice of 
stable perturbation, and it is usually called a {\it $0$-stable invariant} of the original germ $f$. 
We let  $\#\underline{\eta}(f_t)$ denote the number of (multi-)singular points in the  {\em target space}, 
e.g., $\# A_0^3(f_t)$ ($\kappa >0$) means the triple-points in the image of a stable perturbation $f_t$ of $f_0$. 
 
There have been several works for computing $0$-stable invariants via colengths of some algebras under certainly nice situations. 
In general, however, it can be a hard task e.g., in case that the singularity scheme is not Cohen-Macaulay,  that often arises. 

Now we propose a new topological method localizing Thom polynomials 
to compute stable invariants for {\it quasi-homogeneous} map-germs. 
This new approach provides a number of examples with significantly simpler computation where any corank condition is not assumed. 

Let $f: (\C^m, 0) \to (\C^n, 0)$ be a quasi-homogeneous germ 
with  weights $w_1, \cdots, w_m$ and  degrees $d_1, \cdots, d_n\in \Z_{>0}$, i.e., 
there are diagonal representations $\rho_0, \rho_1$ of $T=\C^*$ in the source and the target space, respectively, such that $f=\rho_1(t) \circ f \circ \rho_0^{-1}(t)\; (t \in T)$. 
Suppose that $f$ is $\A$-finitely determined. 
Then its $\A_e$-versal unfolding $F: (\C^{m+k},0) \to (\C^{n+k},0)$ 
is also quasi-homogeneous, so let $r_1, \cdots , r_k$ be the weights of unfolding parameters. 
We introduce vector bundles over $BT=\Proj^\infty$
$$E_0:=\oplus_{i=1}^m \Ost(w_i), \quad 
E_1 :=\oplus_{j=1}^n \Ost(d_j), 
\quad 
E'=\oplus_{i=1}^k \Ost(r_i),$$ 
and we have maps between total spaces of these bundles: 
$$\xymatrix{
&E_0 \ar[r]^f \ar[d]^{i_0}&  E_1 \ar[d]^{i_1}& \\
\underline{\eta}(F) \ar@<-0.3ex>@{^{(}->}[r]  & E_0\oplus E' \ar[r]_F & E_1\oplus E' &\ar@{_{(}->}@<0.3ex>[l] F(\underline{\eta}(F))
}
$$
Given a singularity type $\underline{\eta}$, we apply the Thom polynomial $tp(\underline{\eta})$ for $F$ \cite{Kaz03b, Ohmoto24} and compute the intersection of $tp(\underline{\eta})(c(F))$ with the inclusion map $i_0$ or $i_1$ in $H^*(BT)=\Z[a]$, that leads to the following theorem. 

\begin{thm}\label{counting_0_stable_inv}  \cite{Ohmoto16}
Let $m \le n$ and 
let $f: (\C^m,0) \to (\C^n,0)$  be 
an $\A$-finitely determined quasi-homogeneous map-germ 
with weight $w_i$ and degree $d_j$. 
Given a stable mono/multi-singularity type $\underline{\eta}$ with the expected codimension $n$ in the target, 
the $0$-stable invariant of $f_0$ is computed by 
$$\# \underline{\eta}(f_t) 
= \frac{tp(\underline{\eta})}{{ \deg_1 \underline{\eta}} \cdot w_1\cdots w_m}
= \frac{tp_{\rm target}(\underline{\eta})}{ d_1\cdots d_n}$$
where numerators stand for the coefficient of $a^m$ and $a^n$ 
of the source and target Thom polynomials applied to the universal map 
$f_0: E_0 \to E_1$, respectively.  
\end{thm}

In Example \ref{icis} in \S \ref{higher_tp}, we have seen that the Greuel-Hamm formula of the Milnor number for quasi-homogeneous $\K$-finite map-germs with $\kappa \le 0$ (i.e., ICIS) can be interpreted as a torus localization of the higher Thom polynomial $tp^{\SM}(\mu)$. 
As a counterpart in case of $\kappa >0$,  let us consider $\A$-finitely determined map-germs $f: (\C^n, 0) \to (\C^{n+1},0)$ where the pair $(n, n+1)$ is in Mather's nice range. 
The image of $f$ is a hypersurface with non-isolated singularities -- 
it is known that the homotopy type of the image of any stably perturbation $f_t$ is a bouquet of middle dimensional real spheres, whose number is 
called the {\em image Milnor number}: 
$$\mu_I(f):=(-1)^m(\chi(\Image(f_t))-1).$$ 
The so-called {\em Mond conjecture} says that the $\A_e$-codimension of $f$ $\le$ $\mu_I(f)$, with equality if $f$ is quasi-homogeneous \cite{MondNuno}. 
A crucial difficulty lies in how to calculate  $\mu_I(f)$. 
Our method provides a lot of new examples, i.e., we can compute $\mu_I$ for $\A$-finite quasi-homogeneous germs in low dimensions without any corank condition, by using the same type theorem as above by replacing $tp$ by higher Thom polynomial $tp^{\SM}$ for multi-singularities,  see \cite{Ohmoto16} for the detail in case $n=3$. 
Furthermore, combining this method with some results on several invariants for $\A$-equivalence, 
Pallar\'es Torres and Pe\~nafort Sanchis \cite{PPe21} succeeded to obtain formulas for $\mu_I(f)$ in more cases, e.g., $n=4$ and $5$. 
Their results have a significant impact on the study of the Mond conjecture  \cite[\S 8.8, \S 8.9]{MondNuno}.

\subsection{Positivity and complex hyperbolic geometry}\label{hyperbolic} 
Recently a surprising connection was discovered by B\'erczi \cite{Berczi12, Berczi19a, Berczi19b, Berczi20, Berczi23} between the positivity question of coefficients of Thom polynomials of $A_\mu$-singularities and higher dimensional versions of the Nevanlinna theory in complex analysis. It was quite unexpected, as there had been no interaction between these two subjects before. The discovery was made possible by a new geometric invariant theory for {\em non-reductive groups}, established by B\'erczi-Kirwan \cite{BercziKirwan24}. 

In complex analytic geometry, one of fascinating themes is the {\em hyperbolicity} of a complex projective variety $X$ (in the sense of Brody), i.e., any holomorphic curve from $\C$ to $X$ is constant -- the famous Kobayashi conjecture says that a generic hypersurface $X \subset \Proj^{n+1}$ is hyperbolic if $\deg(X)$ is sufficiently large. It is known that a generic hypersurface $X$ is of general type if  $\deg(X) \ge n+3$. As a generalization, the Green-Griffiths-Lang (GGL) conjecture claims that for every complex projective algebraic variety $X$ of general type there exists a proper algebraic subvariety of $X$ containing all nonconstant entire holomorphic curves $\gamma: \C \to X$. To attack these conjectures, a key construction is a smooth compactification of the bundle of invariant $k$-jet differentials, called the {\em Demailly-Semple tower}, so that global sections of certain line bundles over it give algebraic differential equations which the curves $\gamma$ should satisfy. Then the problem is reduced to the existence of global sections, thus intersection theory comes up here. By this strategy, there have been several attempts to find effective degree bounds for the conjectures in past decades. 

A technical difficulty lies on the fact that the group $J^k\Aut(\C,0)$ of $k$-jets of isomoprhism-germs of $(\C, 0)$ is not reductive, thus 
the geometric invariant theory (GIT) of Mumford can not be applied directly. Recently, B\'erczi and Kirwan developed {\em non-reductive} GIT, and as an application, they obtained a different smooth compactification of the bundle of invariant jet differentials using the algebraic model of B\'erczi-Szenes \cite{BercziSzenes12} for computing the Thom polynomials of $\A_\mu$-singularity types (Theorem \ref{thm_BS}). Namely, the Demailly-Semple tower is replaced by a new fibrewise compactification $\pi: \tilde{X}_k \to X$ of invariant jet differentials, endowed with a tautological line bundle $\tau$ on $\tilde{X}_k$, that leads to the discovery/refinement of  an iterated residue formula for the tautological integrals. 

\begin{thm}[B\'erczi \cite{Berczi12, Berczi19a, Berczi19b, Berczi20}]
\label{GGL1} 
Let $X \subset \Proj^{n+1}$ be a smooth projective hypersurface and let  $u=c_1(\tau)$ and $h = c_1(\pi^*\Ost_X(1))$ denote the first Chern classes of the tautological line bundle on $\tilde{X}_\mu$ and the divisor class on $X$,  respectively. For any homogeneous $\mu$ polynomial $P=P(u,h)$ of degree $\deg(P)=\dim \tilde{X}_\mu =n+\mu(n-1)$, we have 
$$
\int_{\tilde{X}_\mu} P=
\int_X \underset{{\bm{t}}=\infty}{\Res}\;
\frac{\prod_{m<l} (t_m-t_l) Q_k(t_1, \cdots , t_\mu)  P(t_1+ \cdots + t_\mu, h)}
{\prod_{m+r\le l \le \mu} (t_m+t_r-t_l)(t_1\cdots t_\mu)^n}
\prod_{i=1}^\mu {\rm S}\hskip-2pt \left(\frac{1}{t_i}\right) dt_i
$$
where ${\rm S}(t^{-1})=1+s_1(X)t^{-1}+s_2(X)t^{-2}+\cdots + s_n(X)t^{-n}$ is the total Segre class of $X$ 
and $Q_\mu$ is the polynomial associated to $A_\mu$-singularity type in Theorem \ref{thm_BS}. 
\end{thm}

This theorem is very effective for the so-called {\em polynomial GGL conjecture}. For instance, in \cite{Berczi19a, Berczi19b} B\'erczi proved the polynomial GGL conjecture modulo the {\em positivity conjecture} of R. Rim\'anyi,  posed at the very beginning \cite{Rimanyi01}, which claims that there is no negative coefficients in the Chern monomial basis expansion of $tp(A_\mu)$ (this is a different version of positivity from the one in Remark \ref{positivity}). 
After improvements by combining it with recent results in complex hyperbolic geometry,  finally the following theorems have been obtained: 

\begin{thm}[B\'erczi-Kirwan \cite{BercziKirwan24}]
{\bf (Polynomial Green-Griffiths-Lang theorem for projective hypersurfaces)} Let $X \subset \Proj^{n+1}$ be a generic smooth projective hypersurface of  $\deg(X) \ge 16n^3(5n + 4)$. Then there is a proper algebraic subvariety $Y \subset X$ containing all nonconstant entire holomorphic curves in X.
\end{thm}

\begin{thm}[B\'erczi-Kirwan \cite{BercziKirwan24}]
{\bf (Polynomial Kobayashi theorem)} 
A generic smooth projective hypersurface $X \subset \Proj^{n+1}$ of $\deg(X) \ge 16(2n-1)^3(10n-1)$ is Brody hyperbolic.
\end{thm}

These results have made considerable impact in both (mutually far) communities of singularity theory and complex geometry. 
The connecting point between two subjects is intersection theory on Hilbert schemes and non-reductive GIT. 
That also indicates a new future direction for researches in Thom polynomial theory.

\section{Real singularities}\label{real_sing}
We have discussed Thom polynomials in complex analytic/algebraic category, but it makes sense for real $C^\infty$-category as well -- e.g., the desingularization method in \S \ref{desing} commonly works in both categories \cite{AGLV}. There are many related topics in differential topology, but here we focus mainly on a general statement of the existence of real Thom polynomials $tp^\R(\eta)$, i.e., the counterpart to Theorem \ref{tp} in the setting of real singularities. Actually, real singularities have their own particular features caused by a finer classification of $C^\infty$ map-germs, real (semi-)algebraic/analytic properties,  (co-)orientation problem and complex conjugacy. 

\subsection{Vassiliev complex and real Thom polynomials}\label{vassiliev}
The classification machinery (with respect to $\K$ and $\A$-equivalence) is the same as the complex case (\S \ref{classification}), but the real case is more subtle and the classification list becomes larger than the complex case; for instance, there are two different $\K$-orbits with $\kappa=0$, 
$I_{2,2}: (x^2+y^2, xy)$ and $\II_{2,2}: (x^2-y^2, xy)$ as real map-germs, 
but they are the same in complex case. 

Our first task is to define the singularity loci class of prescribed type in (Borel-Moore) homology for real $C^\infty$-maps, that is not straightforward at all. 
Put $V:=J^k(m,n)$ and $G:=J^k\mathcal{G}$, the real algebraic group of $k$-jets of elements in the equivalence group in real case. 
Note that $G$ is homotopy equivalent to the product of orthgonal groups $O(m) \times O(n)$. Since the $G$-action on $V$ is real algebraic, every $G$-orbit $\eta$ (or modulus of $G$-orbits) in $V$ is a locally-closed semi-algebraic subset, and its closure $\overline{\eta}$ in fine topology is also semi-algebraic in general. A semi-algebraic set may not have the $\Z_2$-fundamental cycle in the Borel-Moore homology, while any real algebraic set does  \cite[Prop.~11.3.1]{BCR}.  Also we may ask for the fundamental cycle with $\Z$-coefficients, that is related to the existence of $G$-invariant orientation of the normal bundle of $\eta$ (called the {\em $G$-coorientation} of $\eta$). In this case, we may replace $\mathcal{G}$ by the one preserving the orientations of the source and the target space, then $G$ is  homotopy equivalent to the product of rotation groups $SO(m) \times SO(n)$. For instance, each orbit of $I_{2,2}$ and $\II_{2,2}$ ($\kappa=0$) does not have the fundamental cycle, while the closure of the union of the two orbits forms a real algebraic set, so it has the $\Z_2$-fundamental cycle. Moreover, both orbits are $G$-coorientable and the union admits a $\Z$-fundamental cycle. 

Let us introduce a general framework. 
Suppose that we are given a $G$-invariant semi-algebraic Whitney stratification of $V$. Let $\Gamma_p$ denote the union of strata of codimension $\ge p$ and put $V_p:=V \smallsetminus \Gamma_{p+1}$. 
Then we have a filtration of invariant open sets 
$$V_0 \subset V_1 \subset V_2 \subset \cdots \subset V.$$
It yields the $G$-equivariant spectral sequence $E_*^{*, *}$, which converges to $H_G^*(V)=H^*(BG)$ with the coefficient ring $R$ (e.g., $\Z_2$ or $\Z$).  
The $E_1$-term is 
$$E_1^{p, q}:= H_G^{p+q}(V_p, V_{p-1}; R) \simeq H_G^q(\Gamma_p \smallsetminus \Gamma_{p+1}; R).$$  
In particular, $E_1^{p, 0}$ is freely generated by the strata of codimension $p$ (and $G$-coorientable if $R=\Z$), and 
$\delta: E_1^{p, 0} \to E_1^{p+1, 0}$ is the connecting homomorphism of the long exact sequence for the triple $(V_{p+1}, V_p, V_{p-1})$, which is determined by analyzing the adjacency relation among strata. The cochain complex $(E_1^{*, 0}, \delta)$ is called the {\em Vassiliev complex} (with coefficients in $R$) associated to the stratification \cite{Vassiliev, AGLV, Kaz98}. Intuitively, a $p$-dimensional cohomology class $\omega$ of this cochain complex means a geometric cycle in $V$ defined by a linear combination of strata of codimension $p$, say $\omega=[\sum n_i \eta_i]$ ($n_i \in R$),  and thus it defines an element of $H_G^p(V)=H^p(BG; R)$. This element should be the `real Thom polynomial' associated to the class $\omega$ for $C^\infty$-maps. 
Originally, such an abstract complex was introduced and computed for the classification of function-germs by Vassiliev in his study on Lagrange and Legendre Thom polynomials \cite{Vassiliev, AGLV}, see \S \ref{legendre}. 
For $\K$-classification of $C^\infty$-map-germs with $\kappa=0$ up to the codimension $8$,  the $\Z_2$-Vassiliev complex has been determined in Ohmoto \cite{Ohmoto94}, and the $\Z$-Vassiliev complex, which is much more difficult, has been computed in Feh\'er-Rim\'anyi \cite{FR02} using the interpolation method partly. 

For $\K$-classification, the stabilization property (Lemma \ref{stabilization}, as $m \to \infty$ with $\kappa$ fixed) also holds in real case, that leads to the following theorem. Note that 
$$\textstyle H^*(BO; \Z_2)=\Z_2[w_1, w_2, \cdots], \quad H^*(BSO; \Z[\frac{1}{2}])=\Z[\frac{1}{2}][p_1, p_2,  \cdots]$$ 
where $w_i$ and $p_i$ are the Stiefel-Whitney class and the Pontrjagin class, respectively.

\begin{thm} \label{tp2} 
\cite{Thom55, Damon, Ohmoto94, Kaz98, FR02}. 
Fix $\kappa \in \Z$ and consider $\K$-classification of $C^\infty$-map-germs $(\R^*, 0) \to (\R^{*+\kappa}, 0)$. Let $\omega=[\sum  \eta_i]$ be a $p$-dimensional cohomology class for the $\Z_2$-Vassiliev complex. Then there exists a unique polynomial $tp^\R_2(\omega) \in \Z_2[w_1, w_2, \cdots]$ such that for any appropriately generic $C^\infty$-map $f: M^m \to N^{m+\kappa}$, 
the closure of the union of $\eta_i(f)$, denoted by $S_\omega(f)$, forms a $\Z_2$-geometric cycle in $M$ and 
is expressed by 
$${\rm Dual}_2 [S_\omega(f)] = tp^\R_2(\eta)(w(f)) \;\;   \in H^p(M; \Z_2)$$
evaluated by the Stiefel-Whitney class $w_i=w_i(f)=w_i(f^*TN-TM)$. 
Also for a cohomology class $\omega$ for the $\Z$-Vassiliev complex, a similar statement holds: 
for any generic map $f: M^m \to N^{m+\kappa}$ between orientable manifolds, the locus $S_\omega(f)$ forms a $\Z$-cycle in $M$ and is expressed by a universal polynomial $tp^\R(\omega)$ in Pontrjagin classes $p_i=p_i(f)=p_i(f^*TN-TM)$ modulo $2$-torsion elements. 
\end{thm}

We remark that the interpolation method does not work effectively for computing real mod.2 Thom polynomials $tp^\R_2$, because the restriction equations mod.2 often become trivial (but see Remark \ref{BH} below). Computing real integer Thom polynomials  $tp^\R$ in $p_i$'s is a more difficult problem. 

\begin{rem}\upshape \label{BH} 
{\bf (Borel-Haefliger theorem)} 
As mentioned above,  the union of $\K$-orbits of type $I_{2,2}$ and $\II_{2,2}$ with $\kappa=0$ forms a cycle of the $\Z$-Vassiliev complex; actually its closure in $V$ coincides with the Thom-Boardman stratum $\overline{\Sigma^2}$. Hence, from results by desingularization (which works in both real and complex cases)  \cite{Porteous, Ronga72},  the integer Thom polynomial in the orientable setting is given as 
$$tp^\R(I_{2,2}+\II_{2,2})=tp^\R(\Sigma^2)=p_1$$ 
while forgetting the orientation,  the mod.2 Thom polynomial is 
$$tp^\R_2(I_{2,2}+\II_{2,2})=tp^\R_2(\Sigma^2)=w_2^2+w_1w_3$$ 
(this is consistent; if $M, N$ are orientable, then $w_1=w_1(f^*TN-TM)=0$ and the mod.2 reduction of $p_1$ gives $w_2^2$).  
Also the Thom polynomial for complex singularity type $I_{2,2}$ with $\kappa=0$ is 
$$tp(I_{2,2})=tp(\Sigma^2)=c_2^2+c_1c_3.$$
Note that replacing $c_i$ in $tp(I_{2,2})$ by $w_i$, we recover $tp^\R_2(I_{2,2}+\II_{2,2})$. 
This is generically true. Namely, if a $\K$-orbit $\eta_\C$ of complex map-germs admits several real forms $\eta_i$ such that their complexification is $\eta_\C$ (including the case of a single real form $\eta$), then it defines a $\Z_2$-Vassiliev cycle $\omega=[\sum \eta_i]$, 
and the mod.2 Thom polynomial $tp^\R_2(\omega)$ is obtained from $tp(\eta_\C)$ by replacing $c_i \mapsto w_i$ and by taking coefficients mod.2 (Borel-Haefliger  \cite{BorelHaefliger} and also Matszangosz \cite{Akos}).  
\end{rem}

\subsection{Global topology of $C^\infty$ maps -- Miscellaneous} 
A variant of Vassilev complex for {\em $\A$-classification of multi-singularities of maps} has been introduced for finding topological invariants for stable maps in low dimensions, called {\em local first-order Vassilev-type invariants}, e.g., Arnold's $J^\pm$ invariants for plane immersed curves \cite{Arnold}, Goryunov's invariants for stable maps from a surface into $\R^3$ \cite{Goryunov}, and  many others (see \cite{AicardiOhmoto, Ohmoto12c} and references therein). Putting finer semi-global data for mappings under consideration, an {\em enriched version} of the Vassiliev complex has been considered in Saeki's studies on topology of $4$-manifolds via stable maps into low dimensional space \cite{Saeki, SaekiYamamoto}. That is also related to Gromov's complexity of $C^\infty$-maps. 
One of central problems in global topology of $C^\infty$ maps is elimination of singularities of prescribed type -- historically it arose in connection with the Poincar\'e conjecture. In this context, the h-principle  and cobordism of singular maps take important roles, where the real Thom polynomial of $\eta$ is regarded as the primary obstruction for removing $\eta$-singular points of generic $C^\infty$ maps $f: M \to N$ via homotopy or surgery, see e.g.,  \cite{AGLV, RimanyiSzucs, Szucs, SSS}. 
Also  `relative Thom polynomials' for singular Seifert maps have long been considered implicitly (cf. \cite{Tanabe}); e.g., Takase  \cite{Takase} has studied oriented bordisms of immersions of $7$-manifolds to $\R^8$ in relation to Thom polynomials $tp^\R$ with integer coefficients given in \cite{FR02}. 
For an immersion from a $C^\infty$ manifold $L$ into an almost complex manifold $M$, there are enumerative formulas for complex tangent points or $\R\C$ singularities  \cite{Lai}, and elimination problem for complex tangents has recently been discussed \cite{Elgindi, KasuyaTakase}. 
For a singular Lagrange immersion of $L$ into a symplectic manifold $M$ equipped with almost complex structure, there are Thom polynomials for open Whitney umbrellas \cite{Ohmoto96}. 
Those singularities and complex tangents are characterized by (higher) intrinsic derivatives of the $C^\infty$-bundle map $TL\otimes \C \to TM$, and hence the  Thom polynomials may be written in terms of $c_i(f^*TM-TL\otimes \C)$. 
Higher mod.2 Thom polynomials in real case (i.e., real $Z_2$-counterpart to SSM Thom polynomials) should be based on D. Sullivan's Stiefel-Whitney classes for real algebraic varieties \cite{Sullivan}. There had been only a few related works for real Morin maps, such as Fukuda \cite{Fukuda} and Nakai \cite{Nakai}, but it has potential for further researches, e.g.,  a recent work of Feh\'er-Matszangosz \cite{FM25}. 

There is still a lot of room for further researches on real Thom polynomials (especially in $K$-theory, cobordism theory, Chow-Witt groups, and so on), although the story has a long history since 1950s. Unfortunately several interesting topics are dropped and references cited here may be inadequate, thus readers are advised to directly refer to appropriate literature.

\subsection*{Acknowledgement} 
The author is very grateful to the editors of this series, and to the referee for careful reading and comments for improving the draft. 
On this occasion, the author personally wants to express his deepest gratitude to his teacher Takuo Fukuda, who made a profound study on Isotopy Lemma under Ren\'e Thom in 70's, for guiding him to singularity theory in his student time. Also, the author thanks many friends and colleagues, especially  Laszlo Feh\'er, Maxim Kazarian, Rich\'ard Rim\'anyi, Andrzej Weber for helpful discussions on Thom polynomials for a long time over 20 years. 
This work was partly supported by JSPS KAKENHI Grant Numbers JP17H02838 and JP23H01075.



\end{document}